\documentclass[11pt, reqno]{amsart}
\usepackage{amsmath,amsthm,amssymb}

\theoremstyle{plain}
\newtheorem{thm}{Theorem}[section]
\newtheorem{prop}[thm]{Proposition}
\newtheorem{lem}[thm]{Lemma}
\newtheorem{cor}[thm]{Corollary}

\theoremstyle{definition}

\newtheorem{rem}[thm]{Remark}
\newtheorem{defn}[thm]{Definition}
\newtheorem{eg}[thm]{Example}
\newtheorem{subtitle}[thm]{}
\newtheorem{ex}{Exercise}[section]
\numberwithin{equation}{section}

\def\a{\alpha}
\def\b{\beta}
\def\d{\delta}

\def\e{\epsilon}
\def\g{\gamma}

\def\K{\nabla}
\def\l{\lambda}

\def\n{\,\vert\,}
\def\o{\theta}
\def\w{\omega}
\def\W{\Omega}

\def\cg{{\mathcal{G}}}

\def\co{{\mathcal{O}}}

\def\li{\langle}
\def\ri{\rangle}
\def\n{\ \vert\ }
\def\tr{{\rm tr}}
\def\bs{\bigskip}
\def\ms{\medskip}
\def\ss{\smallskip}

\def\ni{\noindent}
\def\ti{\tilde}
\def\p{\partial}

\def\Im{{\rm Im\/}}
\def\I{{\rm I\/}}
\def\II{{\rm II\/}}
\def\diag{{\rm diag}}
\def\ad{{\rm ad}}

\def\R{\mathbb{R} }
\def\C{\mathbb{C}}
\def\H{\mathbb{H}}

\def\fg{\mathfrak{G}}

\newcommand{\beq}{\begin{equation}}
\newcommand{\eeq}{\end{equation}}
\newcommand{\beg}{\begin{eg}}
\newcommand{\eeg}{\end{eg}}
\newcommand{\bthm}{\begin{thm}}
\newcommand{\ethm}{\end{thm}}
\newcommand{\bprop}{\begin{prop}}
\newcommand{\eprop}{\end{prop}}
\newcommand{\bcor}{\begin{cor}}
\newcommand{\ecor}{\end{cor}}
\newcommand{\blem}{\begin{lem}}
\newcommand{\elem}{\end{lem}}
\newcommand{\bca}{\begin{cases}}
\newcommand{\eca}{\end{cases}}
\newcommand{\brem}{\begin{rem}}
\newcommand{\erem}{\end{rem}}
\newcommand{\bpm}{\begin{pmatrix}}
\newcommand{\epm}{\end{pmatrix}}
\newcommand{\bbm}{\begin{bmatrix}}
\newcommand{\ebm}{\end{bmatrix}}
\newcommand{\bvm}{\begin{vmatrix}}
\newcommand{\evm}{\end{vmatrix}}
\newcommand{\bdefn}{\begin{defn}}
\newcommand{\edefn}{\end{defn}}
\newcommand{\bsub}{\begin{subtitle}}
\newcommand{\esub}{\end{subtitle}}
\newcommand{\bex}{\begin{ex}}
\newcommand{\eex}{\end{ex}}
\newcommand{\ben}{\begin{enumerate}}
\newcommand{\een}{\end{enumerate}}

\def\rd{{\rm \, d\/}}
\def\bu{$\bullet$}

\def\calA{\mathcal A}

\def\calG{\mathcal G}
\def\calJ{\mathcal J}
\def\calK{\mathcal K}
\def\calL{\mathcal L}
\def\calM{\mathcal M}

\def\calP{\mathcal P}
\def\calR{\mathcal R}
\def\calS{\mathcal S}
\def\calU{\mathcal U}
\def\calV{\mathcal V}

\def\fg{\mathfrak{g}}

\def\im{{\bf i}}
\def\ba{{\bf a}}
\def\bc{{\bf c}}

\def\onn{\frac{O(n,n)}{O(n)\times O(n)}}
\def\sdp{\ltimes}

\begin{document}

\title[Geometric Transformations and Soliton Equations]
{Geometric Transformations and Soliton Equations}
\author{Chuu-Lian Terng$^*$}\thanks{$^*$Research supported
in part by NSF Grant DMS-0707132}
\address{Department of Mathematics\\
University of California at Irvine, Irvine, CA 92697-3875}
\email{cterng@math.uci.edu}

\maketitle

\begin{abstract}\vskip 3mm

\noindent
We give a survey of the following six closely related topics: 
(i) a general method for constructing a soliton hierarchy from a splitting of a loop algebra into positive and  negative subalgebras, together with a sequence of commuting positive elements, 
(ii) a method---based on (i)---for constructing soliton hierarchies from a symmetric space, 
(iii) the dressing action of the negative loop subgroup on the space of solutions of the related soliton equation,
(iv) classical B\"acklund, Christoffel, Lie, and Ribaucour transformations for surfaces in three-space and their relation to dressing actions, (v) methods for constructing a Lax pair for the Gauss-Codazzi Equation of certain submanifolds that admit Lie transforms, (vi) how soliton theory can be used to generalize classical soliton surfaces to submanifolds of higher dimension and co-dimension.

\vskip 4.5mm


\end{abstract}

\tableofcontents

\section{Introduction}

Although it is difficult to give a formal definition of soliton equations,  it is generally agreed that a soliton equation is a non-linear wave equation having the following properties (cf. \cite{DS84, AC91, Pa97, TU98}):

\ss\ni {\bf Existence of explicit $n$-soliton solutions} 

A {\it solitary wave\/} is a traveling wave of the form
$u(x,t)=f(x-ct)$ for some smooth function $f$ that decays rapidly as $|x|\to \infty$.  An
$n$-{\it soliton\/} solution is a solution that is asymptotic to  a
nontrivial sum of
$n$ solitary waves
$\sum_{i=1}^n f_i(x-c_it)$ as $t\to -\infty$ and to the sum of the same waves
$\sum_{i=1}^n f_i(x-c_it+r_i)$ with some nonzero phase shifts $r_i$ as
$t\to
\infty$. In other words, during nonlinear interaction, the individual
solitary waves pass through each other, keeping their velocities and
shapes, but with  phase shifts. 

\ss\ni  {\bf ODE B\"acklund transformation}

An ODE B\"acklund transformation is a system of compatible ODEs associated to a given solution of the soliton equation such that solutions of the ODE system are again solutions of the soliton equation. If we apply these transformations to the vacuum solutions repeatedly, then we get explicit multi-soliton solutions. 

\ss\ni  {\bf Bi-Hamiltonian structure and commuting flows}

A pair of Poisson structures $(\{\, , \}_0, \{\, , \}_1)$ on $M$ is called a {\it bi-Hamiltonian structure\/} if $c_0\{\, , \}_0+ c_1\{\, , \}_1$ is a Poisson structure for all constants $c_0, c_1$. A soliton equation is an evolution equation on a function space. One important property is that this function space admits a bi-Hamiltonian structure such that the equation is Hamiltonian with respect to both Poisson structures.  Moreover, one can use these two Poisson structures to construct a hierarchy of commuting Hamiltonian PDEs.  

\ss\ni {\bf Lax pair and inverse scattering}

 A PDE for $q:\R^n\to V$ is said to have a $\calG$-valued {\it Lax pair\/} or a {\it zero curvature formulation\/} if there is a family of $\calG$-valued connection $1$-forms $\o_\l$ on $\R^n$ written in terms of  $q$ and derivatives of $q$ for $\l$ lies in an open subset $\co$ of $\C$ such that the PDE for $q$ is given by the condition that $\o_\l$ is flat for all $\l\in \co$, where $\calG$ is a finite dimensional Lie algebra.  
The Lax pair gives a linear system with a ``spectral parameter'' $\l$. 
The scattering data of a solution is the ``singularity'' of parallel frames of $\o_\l$. The inverse scattering reconstructs the solution from its scattering data (cf. \cite{BC84, TU98}).  

\ss
The above properties will be discussed in more detail in later sections.  
Soliton equations also have algebraic geometric solutions via the spectral curve formulation (cf. \cite{Kri77}), a tau function and a Virasoro action (cf. \cite{Wil91, vM94}). 

\ss\ni {\bf Model soliton equations}

Below are some soliton equations found in 1960s and 70s:
   The Korteweg-de Vries equation (KdV) 
 $$q_t=\frac{1}{4}(q_{xxx}+ 6qq_x),$$
  the non-linear Schr\"odinger equation (NLS) \cite{ZS72}
 $$q_t=\frac{\im}{2}(q_{xx}+ 2|q|^2 q),$$ the modified KdV (mKdV)  
 $$q_t= \frac{1}{4}(q_{xxx} + 6q^2 q_x),$$ the sine-Gordon equation (SGE)
 $$q_{xt}= \sin q,$$ and the $3$-wave equation \cite{ZS79} for $u=(u_{ij})\in su(3)$ with $u_{ii}=0$ for $1\leq i\leq 3$:
 $$(u_{ij})_t= \frac{b_i-b_j}{a_i-a_j} (u_{ij})_x + \frac{b_k-b_j}{a_k-a_j} u_{ik} u_{kj}, \qquad 1\leq i, j, k\leq 3 \,\, {\rm distinct,}$$ where $a_1, a_2, a_3$ are fixed distinct real numbers and $b_1, b_2, b_3$ are fixed real constants.  Although KdV and SGE as soliton equations were discovered in the 1960s and 1970s, they were already studied in the nineteen century.  
 
 \ss\ni {\bf Construction of soliton hierarchy from splittings of Lie algebras}
 
 Zakharov-Shabat found a $sl(2)$-valued Lax pair for NLS in \cite{ZS72},
  Ablowitz-Kaup-Newell-Segur \cite{AKNS74} found $sl(2)$-valued Lax pairs for KdV, mKdV, and SGE, Zakharov-Shabat \cite{ZS79} considered equations admitting a zero curvature formulation depending rationally on $\l$, Adler \cite{Adl79} derived KdV from a splitting of the Lie algebra of pseudo-differential operators on the real line, Kupershmidt-Wilson \cite{KW81} found a $n\times n$ generalization of mKdV,  Drinfeld-Sokolov  \cite{DS84} and Wilson \cite{Wil91}  constructed soliton hierarchies from splitting of loop algebras. These works led to a general method to construct soliton equations from a splitting of Lie algebras.  Many properties of soliton equations can be derived in a unifying way from Lie algebra splittings (cf.  \cite{DS84, Wil91, TU09a}).  
  
  \ss\ni{\bf Soliton hierarchy associated to symmetric spaces}
 
 Given a symmetric space $\frac{U}{K}$, there is a natural Lie subalgebra $\calL$ of the Lie algebra of loops in $\calU\otimes \C$ and a  splitting of $\calL$, where $\calU$ is the Lie algebra of $U$.  We call the soliton hierarchy constructed from this splitting the $\frac{U}{K}$-hierarchy. For example, the $SU(2)$-hierarchy contains NLS, the $\frac{SU(2)}{SO(2)}$-hierarchy contains the mKdV, and the $SU(3)$-hierarchy contains the $3$-wave equation.    If the rank of $\frac{U}{K}$ is $n$, then the first $n$ flows in the $\frac{U}{K}$-hierarchy are PDEs of first order similar to the $3$-wave equation.  We put these first $n$ flows together to construct the $\frac{U}{K}$-system in \cite{Ter97}.  It turns out that many $\frac{U}{K}$-systems are Gauss-Codazzi equations for special classes of submanifolds admitting geometric transforms.  

\ss\ni {\bf Soliton equations in classical differential geometry}
 
Soliton equations were also found in classical differential geometry.  The SGE arose first
through the theory of surfaces of constant Gauss curvature $K=-1$ in $\R^3$, and
the reduced 3-wave equation can be found in Darboux's work \cite{Da1910} on triply
orthogonal coordinate systems of $\R^3$.  In 1906, da Rios, a student of Levi-Civita,
wrote a master's thesis, in which he modeled the movement of a thin vortex
by the motion of a curve propagating in
$\R^3$ along its binormal with curvature as speed. It was  much later, in 1971, that  Hasimoto 
 showed the equivalence of this system with the NLS.   These equations were
rediscovered  independently  of their geometric history. The main contribution
of the classical geometers lies in their methods for constructing explicit
solutions of these equations from geometric transforms.  
  For example:

\ss\ni {\bf $K=-1$ surfaces in $\R^3$, SGE, and B\"acklund transforms} \cite{Ei62}

There is a Tchebyshef line of curvature coordinate system on surfaces in $\R^3$ with $K=-1$ such that the Gauss-Codazzi equation written in this coordinate system is the SGE.   Given a surface $M$ with $K=-1$ in $\R^3$, there is a one parameter family of new surfaces of curvature $-1$ related to $M$ by B\"acklund transformations (a special type of line congruence, see section \ref{cz}). Moreover, this family of new $K=-1$ surfaces can be constructed from a system of ODEs and infinitely many families of explicit solutions of SGE are constructed. 
 
 \ss\ni {\bf Isothermic surfaces in $\R^3$ and Ribaucour transforms} \cite{Da1899}
 
A surface in $\R^3$ is called {\it isothermic\/} if it is parametrized by a conformal line of curvature coordinate system. The Gauss-Codazzi equation written as a first order system is  a soliton equation. 
Given an isothermic surface $M$ in $\R^3$, there is a family of isothermic surfaces related to $M$ by Ribaucour transforms (a special type of sphere congruence, see section \ref{de}).  Moreover, this family of new isothermic surfaces can be constructed by solving a system of compatible ODEs.  

\ss\ni {\bf Higher dimension generalizations via differential geometry}

In late 1970s, S. S. Chern suggested to Tenenblat and the author that the Gauss-Codazzi Equation of $n$-submanifolds in $\R^{2n-1}$ with negative constant sectional curvature might be a new soliton equation in more than two variables.  We found a good coordinate system to write down the Gauss-Codazzi equations in terms of a map from $\R^n$ to $O(n)$ (the generalized sine-Gordon equation GSGE), constructed B\"acklund transformations, a permutability formula, and explicit mutli-soliton solutions  for GSGE in \cite{TenTer80, Ter80}.  Ablowitz, Beals, and Tenenblat \cite{ABT86} constructed a Lax pair for GSGE and used the inverse scattering method to solve the Cauchy problem for GSGE for small rapidly decaying initial data on a non-characteristic line.  Although GSGE is a PDE in $n$ variables, it is really a system of $n$ commuting determined hyperbolic systems in one space and one time variables.   
Tenenblat generalized  B\"acklund theory to other space forms in \cite{Ten85}. Dajczer and Tojeiro constructed Ribaucour transforms for flat Lagrangian submanifolds in $\C^n$ and $\C P^n$ in \cite{DT95, DT00}.  
It turns out that all these geometric equations arise naturally as $\frac{U}{K}$-systems or twisted $\frac{U}{K}$-systems in soliton theory. 

\ss\ni {\bf $\R$-action and associated family}

One reason why many soliton equations arise in submanifold geometry can be seen from the method of moving frames:  A local orthonormal frame $g=(e_1, \ldots, e_{n+k})$ for a submanifold $M^n$ in $\R^{n+k}$ is called {\it adapted\/} if $e_1,\ldots, e_n$ are tangent to $M$.  
The Gauss-Codazzi equation (GCE) for $M$ is given by the flatness for the Maurer-Cartan form $\o=g^{-1} \rd g$.  Consider a class of $n$-submanifolds in $\R^{n+k}$ satisfying a certain geometric condition.  Suppose
\ben
\item[(a)] we can use this geometric condition to find a ``good'' coordinate system on these submanifolds such that its Maurer-Cartan form $\o$ and hence the GCE has specially ``simple'' form,
\item[(b)]  there is an $\R$-action on solutions of the GCE, and we call an orbit of the induced $\R$-action on this class of submanifolds an {\it associated family\/}. 
\een
Then the induced $\R$-action on the Maurer-Cartan form often gives a Lax pair for the Gauss-Codazzi equation, which is one of the characteristic properties of soliton equations.  Thus we call a class of submanifolds {\it soliton submanifolds\/} if its Gauss-Codazzi equation is a soliton equation.  

\ss\ni {\bf Higher dimension generalization via soliton theory}

Constructions and generalization of geometric transforms for soliton surfaces in $\R^3$ to submanifolds in $\R^n$ are beautiful but mysterious and usually are done case by case.  However, geometric transforms for soliton submanifolds in $\R^n$ can be constructed in a unified way from the action of ``simple'' rational loops on the space of solutions of soliton equations and the permutability formula is then a consequence of the geometric transforms being part of a group action.  
If the Gauss-Codazzi equation of a class of surfaces in $\R^3$ admitting geometric transforms is a soliton equation associated to a rank $2$ symmetric space, then we can often use the same type of symmetric space of higher rank to construct a natural generalization of a class of soliton surfaces in $\R^3$ to higher dimension and co-dimension soliton submanifolds.  For example, the Gauss-Codazzi equation for Christoffel pairs of isothermic surfaces in $\R^3$ is the $\frac{O(4,1)}{O(3)\times O(1,1)}$-system \cite{CGS95}, which led to a natural generalization to $k$ tuples of isothermic $k$-submanifolds in $\R^n$ whose equation is the $\frac{O(n+k-1,1)}{O(n)\times O(k-1,1)}$-system. Moreover, the action of rational loops on this $\frac{U}{K}$-system gives rise to natural generalizations of Ribaucour transforms and permutability formulae for these $k$ tuples of isothermic submanifolds in $\R^n$ (cf. \cite{DonTer08a}).

This article is organized as follows: We set up notations for the moving frame method for submanifolds in section 2, review the classical notion of line congruences and geometric B\"acklund transformations for surfaces in $\R^3$ with $K=-1$ and $n$-submanifolds in $\R^{2n-1}$ with constant sectional curvature $-1$ in section 3, and explain the notions of sphere congruences and Ribaucour transforms for isothermic surfaces in section 4. In section 5 we review Combescure transforms, O surfaces, and $k$-tuples in $\R^n$ and the fact that $k$-tuples in $\R^n$ give a natural generalization of isothermic surface theory to arbitrary dimension and co-dimension isothermic submanifolds.  
  In section 6 we derive the Lax pairs for Gauss-Codazzi equations using the moving frame of the associated family for surfaces in $\R^3$ with $K=-1$, isothermic surfaces, $k$-tuples in $\R^n$, and flat Lagrangian submanifolds in $\C^n$. In section 7 we give a brief discussion of the method of constructing soliton hierarchies from splittings of loop algebras and derive formal inverse scattering, commuting flows, and bi-Hamiltonian structure from the splitting. We give definitions of $\frac{U}{K}$-system, twisted $\frac{U}{K}$-system, and the $-1$ flow on the $\frac{U}{K}$-system and their Lax pairs in section 8.  We review the construction of the action of the group of rational maps $f:S^2=\C\cup \{\infty\}\to U_\C$ such that $f(\infty)=\I$ and $f$ satisfies the $\frac{U}{K}$- reality condition on the space of solutions of the $\frac{U}{K}$-system in section 9.  In the final section, we give the relation between the rational loop group action on the space of solutions of $\frac{U}{K}$-system and  geometric transformations of the corresponding soliton submanifolds.
  
  The author selects only few classes of soliton submanifolds in Euclidean space to explain the relation between various geometric transforms and group actions on solutions of soliton equations.  The reader may find more examples of soliton submanifolds of space forms and symmetric spaces in \cite{Ten85, Bo94,  FP96b, BDPT02, MN06, Bra07a, Bra07b}, soliton surfaces in affine geometry in \cite{Bo99, W06},  and soliton submanifolds of conformal geometry in \cite{DonTer08b, BDPP09}. For the theory of soliton equations, we refer the reader to \cite{DS84, AC91, Pa97} and for the theory of transformations we refer the reader to \cite{For94, Ten98, GHZ99}. We also refer to these references for more complete lists of works related to soliton equations and soliton submanifolds.

\bs
\section{The moving frame method for submanifolds} 

 Let  $f:M^n\to\R^{n+k}$ be an immersion, and $(\, ,)$ the standard inner product on $\R^{n+k}$. The first and second fundamental forms $\I, \II$ and the induced normal connection $\K^\perp$ form a complete set of local invariants and they must satisfy the Gauss-Codazzi equations.  Below we set up notations for the method of moving frames of Cartan and Chern. 

Let $g=(e_1, \ldots, e_{n+k})$ be a local orthonormal frame on $M$ such that $e_1, \ldots, e_n$ are tangent to $M$, and let $w_1,\ldots, w_n$ be the $1$-forms on $M$ dual to $e_1, \ldots, e_n$. Then 
\beq\label{aa}
\rd f= \sum_{i=1}^n w_i e_i.
\eeq
Since $g^tg=\I$, the Maurer-Cartan form
$$w=(w_{AB}):= g^{-1} \rd g$$
 is $o(n+k)$-valued.  In other words, $\rd g= gw$, i.e., 
$$\rd e_B= \sum_{A=1}^{n+k} w_{AB} e_A, \quad {\rm or \,\, equivalently, \, \,\,} w_{AB}= (\rd e_B, e_A). $$
We use the following index conventions:
$$1\leq i, j, k\leq n, \quad n+1\leq \a, \b, \g\leq n+k, \quad 1\leq A, B, C\leq n+k.$$
Then $\I, \II, \K^\perp$ are given by
$$\I= \sum_{i=1}^n w_i^2, \quad
 \II= \sum_{i=1, \a= n+1}^{n, n+k} w_i w_{i\a} e_\a,\quad
 \K^\perp e_\a= (\rd e_\a)^\perp= \sum_{\b} w_{\b\a}e_\b,$$
 where $\xi^\perp$ denotes the projection of $\xi$ onto $\nu(M)$ along $TM$.  The shape operator $A_v$ along a normal vector $v\in \nu(M)_p$ is the self-adjoint operator on $TM_p$ defined by $(\II(u_1, u_2), v)= (A_v(u_1), u_2)$ for all $u_1, u_2\in TM_p$.  
 
\blem\label{fi} {\bf Cartan Lemma}\par

 The Levi-Civita connection $1$-form $(w_{ij})_{1\leq i, j\leq n}$ for $\I= \sum_{i=1}^n w_i^2$ is obtained by solving the {\it structure equation\/}:
\beq\label{ad}
\rd w_i= -\sum_{j=1}^n w_{ij}\wedge w_j, \quad w_{ij}+ w_{ji}=0, \quad 1\leq i, j\leq n.
\eeq
\elem

For example, the Levi-Civita connection $1$-form $(w_{ij})$ for a diagonal metric $\I= \sum_{i=1}^n a_i(x)^2 \rd x_i^2$ is 
\beq\label{ae}
w_{ij}= \frac{(a_i)_{x_j}}{a_j} \rd x_i - \frac{(a_j)_{x_i}}{a_i} \rd x_j.
\eeq

\ss\ni {\bf Gauss-Codazzi equations}

Since $w=g^{-1}\rd g$, $w$ is a flat $o(n+k)$-valued connection $1$-form, i.e., $\rd w= - w\wedge w$. Or equivalently,
\beq\label{ab}
  \rd w_{AB} = -\sum_{C} w_{AC}\wedge w_{CB}, \quad 1\leq A\leq n+k.
\eeq
This gives the Gauss-Codazzi-Ricci equation for $M$:
\begin{gather}
\W_{ij}=\rd w_{ij}+ \sum_k w_{ik}\wedge w_{kj}= \sum_\a w_{i\a}\wedge w_{j\a},\label{afg}\\
 \rd w_{i\a} = -\sum_j w_{ij}\wedge w_{j\a} - \sum_\b w_{i\b}\wedge w_{\b\a},\label{afc}\\
\W^\perp_{\a\b}=\rd w_{\a\b} + \sum_{\g} w_{\a\g}\wedge w_{\g\b} = \sum_i w_{i\a}\wedge w_{i\b}, \label{afr}
\end{gather}
where $\W_{ij}$ and $\W_{\a,\b}^\perp$ are the curvature tensors for $\I$ and for the induced normal connection $\K^\perp$ respectively. 

Write $w_{i\a}= \sum_j h^\a_{ij} w_j$. Then $h_{ij}^\a= h_{ji}^\a$ and the matrix for the shape operator$A_{e_\a}$ is $(h^\a_{ij})$ with respect to the tangent basis $e_1, \ldots, e_n$.  The Ricci equation gives 
$$\W^\perp_{\a, \b} =\sum_i w_{i\a}\wedge w_{i\b} = \sum_{i, k, l} h_{ik}^\a h_{il}^\b w_k\wedge w_l.
$$

\ss\ni {\bf Flat and non-degenerate normal bundle}

 The normal bundle is {\it flat\/} if the normal curvature is zero, i.e.,
$\W_{\a\b}^\perp=0$, or equivalently $ [A_{e_\a}, A_{e_\b}]=0$ for all $n+1\leq \a, \b\leq n+k$.    
So the normal bundle is flat if and only if all shape operators commute.  In this case, for fixed $p\in M$, we can find a common eigenbasis for the shape operators $\{A_v\n v\in \nu(M)_p\}$. 

The normal bundle of an $n$-dimensional submanifold in $\R^{n+k}$ is {\it non-degenerate\/} if for each $p$
the space of shape operators $\{A_v\n v\in \nu(M)_p\}$ has dimension $k$.  

\bthm  {\bf Fundamental Theorem of submanifolds in $\R^N$}\label{di} \cite{PalTer88}\par

Let $M$ be an open subset of $\R^n$, and $\eta$ an orthogonal rank $k$ vector bundle on $M$ with an $O(k)$-connection $\ti \K$.  Let $\fg$ be a Riemannian metric on $M$, and $\xi$ a smooth section of $S^2(T^*M)\otimes \eta$.  
 We construct an $o(n+k)$-valued $1$-form as follows:
\ben
\item Choose $1$-forms $w_1, \ldots, w_n$ such that $\fg= \sum_{i=1}^n w_i^2$.
\item Solve $(w_{ij})_{1\leq i, j\leq n}$ from the structure equation \eqref{ad}.
\item Choose a local orthonormal frame $(s_{n+1}, \ldots, s_{n+k})$ for $\eta$. Write the connection $\K^\perp s_\a= \sum_\b w_{\b\a} s_\b$.
\item Write  $\xi= \sum_{\a, i,j} h^\a_{ij} w_i w_js_\a$ with $h_{ij}^\a= h_{ji}^\a$.  Set $w_{i\a}=-w_{\a i}= \sum_j h^\a_{ij} w_j$.
\een 
If $w:=(w_{AB})_{1\leq A, B\leq n+k}$ is a flat $o(n+k)$-valued connection $1$-form, i.e., $\rd w=-w\wedge w$, then given $x_0\in M$, $p_0\in \R^{n+k}$, and an orthonormal basis $\{v_1, \ldots, v_{n+k}\}$ of $\R^{n+k}$,   the following system of first order PDE for $(f, e_1, \ldots, e_{n+k})$ is solvable and has a unique solution defined in an open subset $\co$ of $x_0$ in $M$:
\beq \label{ag}
\bca \rd f= \sum_i w_i e_i,\\
\rd e_A= \sum_B w_{BA} e_B,\\
f(x_0)= p_0, \quad e_A(0)= v_A.
\eca
\eeq
Moreover, 
\ben
\item[(a)] $f:\co\to \R^{n+k}$ is an immersion with $\I= \fg$ and $ \II=\sum h_{ij}^\a w_i w_j e_\a$, 
\item[(b)] $e_\a(x)\mapsto s_\a(x)$ gives a vector bundle isomorphism from $\nu(M)$ to $\eta$ that preserves the orthogonal structure and maps the induced normal connection $\K^\perp$ to $\ti \K$ and $\II$ of $f$ to $\xi$.  
\een
\ethm

\brem The Fundamental Theorem \ref{di} can be formulated as the flatness of a $\calG$-valued connection $1$-form, where $\calG$ is the Lie algebra of the rigid motion group $G$ of $\R^{n+k}$: First note that $G$ can be embedded in $GL(n+k+1)$ by 
$$\phi_{g, v}(x)= gx+ v \mapsto \bpm g &v\\ 0 &1\epm, \qquad g\in O(n+k), v\in \R^{n+k}.$$
The Lie algebra of the rigid motion group is the subalgebra of $gl(n+k+1)$:
$$\calG= \left\{ \bpm A & v\\ 0 &0\epm\, \bigg| \, A\in o(n+k), v\in \R^{n+k}\right\}.$$
The equation for isometric immersion for given $\I, \II, \K^\perp$ is \eqref{ag}, or equivalently
$$\rd \bpm g & f\\ 0 &1\epm = \tau\bpm g & f\\ 0 &1\epm, \qquad {\rm where}\quad 
\tau= \bpm w_{ij} & w_{i\a} & w_i\\ w_{\a i} & w_{\a\b} &0\\ 0&0&0\epm.$$
  This  system is solvable for any initial data $c_0\in O(n+k)$ and $p_0\in \R^{n+k}$ if and only if $\tau$ is flat. Or equivalently, $w_i, w_{AB}$ satisfy the structure equation \eqref{ad} and the Gauss-Codazzi equation \eqref{ab}. 
\erem

\bs
\section{Line congruences and B\"acklund transforms}  \label{cz}

We review the classical notion of line congruences and geometric B\"acklund transforms for $K=-1$ surfaces in $\R^3$ and for $n$-submanifolds in $\R^{2n-1}$ with constant sectional curvature $-1$ (\cite{Ei62, TenTer80, Ter80}).

A {\it line congruence\/} in $\R^{3}$ is a smooth $2$- parameter family of lines, 
$$\ell(x)=\{c(x)+ t v(x)\n t\in \R\}$$  
defined for $x$ in an open subset $\co$ of $\R^2$. 
A surface $f:\co\to \R^{3}$ is called a {\it focal surface\/} of the line congruence $\ell$ if 
$f(x)\in \ell(x)$ and $\ell(x)$ is tangent to $f$ at $f(x)$ for each $x\in \co$. To find a focal surface is to find a function $t:\co\to \R$ such that $f(x)= c(x) + t(x) v(x)$ is an immersion and $v(x)$ is tangent to $f$ at $f(x)$.  This condition is equivalent to 
$$\det(f_{x_1}, f_{x_2}, v)=0,$$
which is a quadratic equation in $t$.  So generically, there are exactly two focal surfaces for a line congruence. Moreover, the two focal surfaces determine the line congruence.   Hence we call a diffeomorphism $\phi:M\to \ti M$ a {\it line congruence\/} if the line jointing $p$ and $\phi(p)$ is tangent to $M$ and $\ti M$ at $p$ and $\phi(p)$ respectively for all $p\in M$.

\ss\ni{\bf $K=-1$ surfaces in $\R^3$ and the sine-Gordon equation} (cf. \cite{Ei62, PalTer88})\par
 
We can use the Codazzi equation to prove that  if $M$ is a surface in $\R^3$ with $K=-1$, then locally there exists a line of curvature coordinate system $(x_1, x_2)$ such that 
\beq\label{ar}
\I= \cos^2 q\, \rd x_1^2+ \sin^2 q\,\rd x_2^2, \quad \II= 2\sin q\cos q\,(\rd x_1^2- \rd x_2^2)
\eeq
for some smooth function $q$.  We call $(x_1, x_2)$ the {\it Tchebyshef line of curvature coordinate system}.  Note that $q$ is the angle between the asymptotic lines.  Let $w_1= \cos q \,\rd x_1$ and $w_2= \sin q\,\rd x_2$. By \eqref{ae}, $w_{12}= -q_{x_2} \rd x_1 - q_{x_1} \rd x_2$. Use $\II$ to see that $w_{13}= \sin q \rd x_1$ and $w_{23}= -\cos q \rd x_2$.  The Gauss-Codazzi equation is given by the flatness of 
\beq\label{cx}
w=(w_{AB})=\bpm 0 & -q_{x_2} \rd x_1 - q_{x_1} \rd x_2 & \sin q \rd x_1\\ q_{x_2} \rd x_1 + q_{x_1} \rd x_2 & 0 & -\cos q \rd x_2\\ -\sin q \rd x_1 & \cos q \rd x_2 &0\epm,
\eeq
which gives the sine-Gordon equation (SGE) 
\beq \label{ay} q_{x_1x_1}-q_{x_2x_2} = \sin q \cos q.\eeq
 Change to light cone coordinates $s, t$: $$x_1= s+t, \quad x_2= s-t.$$
 The fundamental forms \eqref{ar} become
$$\I= \rd s^2 + 2\cos (2q) \rd s\rd t + \rd t^2, \quad \II= 2\sin (2q) \rd s \rd t.$$
   The SGE in $(s,t)$ coordinate system is
\beq \label{ba} q_{st}= \sin q\cos q.\eeq
We call $(s,t)$ the {\it Tchebyshef asymptotic coordinate system\/}. 

\bdefn {\bf B\"acklund transformation}

 A line congruence  $\phi:M\to M^*$  is called a {\it B\"acklund transformation\/} (BT) with constant $\o$ if  for any $p\in M$, 
 the distance between $p$ and $p^*= \phi(p)$ is $\sin \o$, and the angle between the normal line of $M$ at $p$ and the normal line of $M^*$ at $p^*$ is equal to
$\o$.
\edefn

 \bthm\label{az} {\bf B\"acklund Theorem}\par 

If $\phi:M\to M^*$ is a B\"acklund transformation with constant $\o$, then both $M$ and $M^*$ have constant Gaussian curvature $K=-1$ and $\phi$ preserves Tchebyshef line of curvature and asymptotic coordinates.  
Conversely, given a surface $M$ in $\R^3$ with $K=-1$, a constant $0<\o<\pi$, $p_0\in M$, and $v_0\in TM_{p_0}$ a unit vector, then there exist a unique surface $M^*$ and a B\"acklund transformation $\phi:M\to M^*$ with constant $\o$ such that $\phi(p_0)= p_0+ \sin \o v_0$.  
\ethm

Analytically to find a BT  $\phi$ with constant $\o$ for a given $K=-1$ surface $M$ in Theorem \ref{az} is to find a unit tangent field $v$ on $M$ such that $\phi(x)= x+ \sin\o v(x)$ is a BT.  Let $e_i$ denote the unit principal directions for $i=1, 2$ and write $v= \cos q^* e_1+ \sin q^* e_2$, then the condition that $\phi$ is a BT with constant $\o$ is equivalent to $q^*$ solving a system of compatible first order ODEs: 

\bthm {\bf ODE B\"acklund transform}\hfil\par

 Given $q(s,t)$ and a non-zero real constant $\mu$, the following system is solvable for $q^*$
\beq\label{aj}
{\rm BT}_{q, \mu}\quad \bca
(q^*+ q)_s=\mu\sin (q^*-q), &\cr
(q^*-q)_t= \frac{1}{\mu} \sin(q^*+q),
\eca
\eeq
if and only if $q$ is a solution of the SGE \eqref{ba}. Moreover, if $q$ is a solution of the SGE then a solution  $q^*$ of \eqref{aj} is again a solution of the SGE.  \ethm

The parameter $\o$ for geometric BT in Theorem \ref{az} and constant $\mu$ in system \ref{aj} are related by 
$\mu= \tan \frac{\o}{2}$.

 Given a solution $q$ of SGE, we can solve the system BT$_{q,\mu}$ to get a family of new solutions of SGE.  If we
apply this method again, then we get a second family of solutions.  This gives infinitely many families of solutions  from a given solution of SGE.   
For example,  the constant function $q=0$ is called the {\it trivial or vacuum solution\/} of the SGE.  The system
BT$_{0,\mu}$ is
$$ \bca
\a_s= \mu \sin \a,&\\
\a_t= \frac{1}{\mu} \sin \a.
\eca $$
It has an explicit solution
\beq\label{by}
\a(s,t)= 2\tan^{-1}\left(e^{\mu s + \frac{1}{\mu} \ t}\right).
\eeq
We can solve B\"acklund transformation BT$_{\a,\mu_1}$ to get another family of solutions. However, BT$_{\a,\mu_1}$ is not as easy to solve as BT$_{0,\mu}$.  
But instead of solving BT$_{\a,\mu_1}$ we can use the following Theorem:

\begin{thm} \label{al} {\bf Bianchi Permutability Theorem}\par

Let $0<\o_1,\o_2<\pi$ be constants such  that $\sin^2\o_1\not= \sin^2 \o_2$, and $\ell_i:M_0\to
M_i$ B\"acklund transformations with constant $\o_i$ for $i=1, 2$.  Then there exist a unique
surface $M_3$ and B\"acklund transformations $\ti\ell_1:M_2\to M_3$ and $\ti \ell_2:M_1\to
M_3$ with constant $\o_1, \o_2$ respectively such that $\ti \ell_1\circ \ell_2= \ti \ell_2\circ \ell_1$.  Moreover, if $q_i$ is the
solution of the SGE corresponding to $M_i$ for $0\leq i\leq 3$, then
\beq\label{ak}
\tan\left(\frac{q_3-q_0}{2} \right)= \frac{\mu_1+\mu_2}{\mu_1-\mu_2}
\tan\left(\frac{q_1-q_2}{2}\right),
\eeq
where $\mu_i = \tan \frac{\o_i}{2}$.  
\end{thm}

\ss\ni {\bf Global verses local}

It follows from the Fundamental Theorem of Surfaces in $\R^3$ that there is a bijective correspondence between  solutions $q$ 
of the SGE \eqref{ay} satisfying $\Im(q)\subset (0,\frac{\pi}{2})$ and local surfaces in $\R^3$ with $K=-1$ up to rigid motions. So we can construct infinitely many families of $K=-1$ surfaces in $\R^3$ by solving compatible systems of ODEs. 
Note that 
if $q:\R^2\to \R$ is a smooth solution of SGE such that $\sin q\cos q$ is zero at a point $p_0$, then although the map $f$ constructed from the Fundamental Theorem of Surfaces in $\R^3$ fails to be an immersion at $p_0$, it is smooth at $p_0$, $df_{p_o}$ has rank $1$ and the tangent bundle is smooth at $p_0$.  Thus global solutions of SGE give $K=-1$ surfaces in $\R^3$ with cusp singularities but smooth tangent bundle.  This is a common phenomenon for soliton submanifolds: Although the Cauchy problem for small norm initial data can be solved globally, the corresponding soliton submanifolds often are only defined locally.

\ss\ni {\bf Explicit multi-soliton solutions for the SGE}

Write the solutions $\a$ of BT$_{0, \mu}$ given in \eqref{by} in space-time coordinates $x_1= s+t$ and $x_2=s-t$ to get
$\a(x_1, x_2)= 2\tan^{-1} e^{\csc\o x_1 - \cot \o x_2}$.  So 
$$\a_{x_1} = \frac{2\csc\o e^{\csc\o x_1- \cot\o x_2}}{1+ e^{2(\csc\o x_1- \cot\o x_2)}}.$$
Note that $\a$ is a traveling wave solution and $\a_{x_1}$ decays to zero as $|x_1|\to \infty$.  Hence SGE viewed as an equation of $\a_{x_1}$ has solitary wave solutions.  These are the $1$-soliton solutions of the SGE.  If we apply permutability formulae to these $1$-solutions, then we get $2$-soliton solutions. Moreover,  these solutions are asymptotically equal to a sum of two solitary waves as $x_2\to -\infty$ and to the sum of the same two solitary waves as $x_2\to \infty$ but with phase shifts (cf. \cite{DaiTer05}).  Explicit multi-soliton solutions of SGE can be obtained by applying permutability formulas repeatedly.  

\ms
\ni {\bf Lie or Lorentz transform}

Lie observed that SGE is invariant under the Lorentz transformations, which are called Lie transforms:
If $q(s. t)$ is a solution of SGE \eqref{ba} and $r$ a non-zero real constant, then  $\ti q(s,t):= q(rs, r^{-1}t)$ is also a solution of SGE.

\ms
\ni {\bf Associated family of $K=-1$ surfaces in $\R^3$}

Given a $K=-1$ surface $M$ in $\R^3$, let $q(s,t)$ denote the corresponding solution of the SGE, $\l\in \R$ a non-zero constant, and $q^\l(s,t)= q(\l s, \l^{-1}t)$.  The family of $K=-1$ surfaces in $\R^3$ corresponding to SGE solution $q^\l$ is called the {\it associated family of $K=-1$ surfaces in $\R^3$ containing $M$}. In section \ref{da}, we will use the moving frame of this associated family to derive the standard Lax pair for SGE.

\ms\ni{\bf $n$-submanifolds in $\R^{2n-1}$ with sectional curvature $-1$ and GSGE} 

The hyperbolic $n$-manifold $\H^n$ is the simply connected, complete, $n$-dimensional Riemannian manifold with constant sectional curvature $-1$. 
\'E. Cartan proved that $\H^n$ can not be locally isometrically immersed in $\R^{2n-2}$, but can be locally isometrically immersed in $\R^{2n-1}$ and the normal bundle of such immersions must be flat (\cite{Car19}).  Moore  used Codazzi equations to prove the existence of  line of curvature coordinate systems on such immersions,  a slight improvement of Moore's result was given in \cite{Ter80} to get an analogue of Tchebyshef line of curvature coordinate systems,  and the corresponding Gauss-Codazzi equation is called the {\it generalized sine-Gordon equation\/} (GSGE).  B\"acklund theory was generalized to GSGE in \cite{TenTer80, Ter80}. 

\bthm  \label{bn}
Let $M^n$ be a simply connected submanifold of $\R^{2n-1}$ with constant sectional curvature $-1$. Then the normal bundle $\nu(M)$ is flat and  there exist coordinates $(x_1, \ldots, x_n)$, an $O(n)$-valued map $A=(a_{ij})$, and parallel normal frames $e_{n+1}, \ldots, e_{2n-1}$ such that the first and second fundamental forms are of the form
$$\I= \sum_{i=1}^n a_{1i}^2 dx_i^2, \quad \II= \sum_{i=1, j=2}^n a_{1i} a_{ji} dx_i^2 e_{n+j-1}.$$
We call $x$ the Tchebyshef line of curvature coordinate system for $M$.
\ethm

To write down the Gauss-Codazzi equation for these immersions we  set 
\begin{gather}
w_i=a_{1i} dx_i, \quad 1\leq i\leq n, \label{bo1}\\
w_{i, n+j-1}=-w_{n+j-1, i}= a_{ji}dx_i.\label{bo3}
\end{gather}
By \eqref{ae}, $w_{ij}= f_{ij} \rd x_i - f_{ji} \rd x_j$, where
\beq\label{bo2}
 f_{ij}= \bca\frac{(a_{1i})_{x_j}}{a_{1j}}, & i\not= j,\\ 0,& i=j,\eca.
 \eeq 
Set $F=(f_{ij})$. Then 
\beq\label{bo4}
\w=(w_{ij})_{i, j\leq n}=\d F- F^t\d, \quad \d=\diag(\rd x_1, \ldots, \rd x_n)
\eeq 
is the Levi-Civita $o(n)$-connection of the induced metric $\I$. 
The Gauss-Codazzi equation and the structure equation give
\beq\label{bp1}
\bca
dw+ w\wedge w= - \d A^t e_{11}A\d,\\
(a_{ki})_{x_j}= f_{ij} a_{kj}, & 1\leq i\not=j\leq n, 1\leq k\leq n,
\eca
\eeq
where $e_{11}$ is the $n\times n$ matrix with all entries zero except the $11$-th entry is $1$.
Or equivalently, it is the second order PDE system for the $O(n)$-valued map $A=(a_{ij})$:
\beq\label{bp2}
\bca (f_{ij})_{x_j} + (f_{ji})_{x_i} + \sum_k f_{ik}f_{jk} = a_{1i} a_{1j}, & i\not= j,\\
(f_{ij})_{x_k} = f_{ik}f_{kj}, & i, j, k \, {\rm distinct,}\\
(a_{ki})_{x_j} = a_{kj} f_{ij}, & i\not= j, \, \forall k.\eca
\eeq
This is the {\it GSGE\/}, and when $n=2$, it is the SGE. 

Since $\sum_{i=1}^n a_{ki}^2=1$, 
$$a_{ki} (a_{ki})_{x_i} = -\sum_{j\not= i} a_{kj} (a_{kj})_{x_i} = -\sum_{j\not= i} a_{kj} f_{ji} a_{ki}.$$
So we have 
\beq\label{fd}
(a_{ki})_{x_i}= -\sum_j a_{kj} f_{ji}.
\eeq  
It follows from \eqref{fd} and the third equation of \eqref{bp2} that 
$$\rd A= A( \d F^t - F\d).$$ So \eqref{bp1} is equivalent to 
\beq\label{bp}
\bca dw+ w\wedge w= - \d A^t e_{11}A\d, & {\rm where\,\, } w= \d F- F^t \d,\\
A^{-1}dA= \d F^t - F\d\eca
\eeq 
Note that we associate to an $n$-submanifold of $\R^{2n-1}$ three flat connections: the flat $o(n)$-connection  $\d F^t- F\d$, the flat $o(n,1)$-connection 
$$\bpm \d F- F^t \d & \xi^t\\ \xi &0\epm, \quad {\rm where\,\,} \xi= (w_1, \ldots, w_n),$$ and the flat $o(2n-1)$ Maurer-Cartan form $(w_{AB})_{A, B\leq 2n-1}$.   

To generalize B\"acklund transformations to higher dimensions, we first recall the notion of  $k$ angles between two $k$-dimensional linear subspace $V_1$ and $V_2$ of a $2k$-dimensional inner product space $(V, (\, , ))$: Let $\pi$ denote the orthogonal projection of $V$ onto $V_1$. Define a symmetric bilinear form on $V_2$ by $\li v_1, v_2\ri = (\pi(v_1), \pi(v_2))$. Then there is a self-adjoint operator $A$ on $V_2$ such that $\li v_1, v_2\ri = (A(v_1), v_2)$.  The $k$ angles between $V_1$ and $V_2$ are $\o_1, \ldots, \o_k$ if $\cos^2\o_1, \ldots, \cos^2\o_k$ are the eigenvalues of $A$.  

\bdefn
Let $M, M^*$ be two $n$-dimensional submanifolds of $\R^{2n-1}$ with flat normal bundle.  A diffeomorphism $\phi:M\to M^*$ is called a {\it B\"acklund transformation with constant $\o$\/} if for all $p\in M$
\ben
\item the line joining $p$ and $p^*=\phi(p)$ is tangent to $M$ at $p$ and to $M^*$ at $p^*$,
\item $||\overline{pp^*}||=\sin \o$,
\item the $(n-1)$ angles between the normal space $\nu(M)_p$ and $\nu(M^*)_{p^*}$  are all equal to the constant $\o$ (note that these normal spaces are two $(n-1)$ dimensional linear subspaces of the $(2n-2)$ dimensional subspace of $\R^{2n-1}$ that is perpendicular to  $p-p^*$).
\een
\edefn

Let $\ell(p)$ denote the line in $\R^{2n-1}$ through $p$ and $\phi(p)$ for a B\"acklund transformation $\phi:M\to M^*$. Then condition (1) says that $\ell$ is an $n$-parameter family of lines in $\R^{2n-1}$ (i.e., an $n$-dimension line congruence in $\R^{2n-1}$) and  $M, M^*$ are focal surfaces of $\ell$.  

\bthm 
If $\phi:M\to M^*$ is a B\"acklund transformation for $n$-dimensional submanifolds in $\R^{2n-1}$ with constant $\o$, then both $M, M^*$ have constant sectional curvature $-1$.  Moreover, $\phi$ maps Tchebyshef line of curvature coordinate system of $M$ to that of $M^*$.
\ethm

Let 
$$\I_{k, n-k}= \diag(\e_1, \ldots, \e_n), \quad {\rm where\,\,} \e_i=1 \, {\rm for\, } i\leq k, \e_i=-1 \, {\rm for\, } k< i\leq n.$$

B\"acklund transform analytically gives

\bthm\label{bd}
Given a smooth $A:\R^n\to O(n)$ and real non-zero constant $\l$,  the following system for $X:\R^n\to O(n)$,
\beq \label{em}
{\rm BT}_{A,\l}: \qquad \rd X = X\d A^t D_\l X - X\w -D_\l A\d,  
\eeq 
is solvable if and only if $A$ is a solution of GSGE, where $D_\l=\frac{(\l\I + \l^{-1}J)}{2}$, $J=\I_{1, n-1}$.  Moreover, the solution $X$ is again a solution of GSGE.  
\ethm

The constant $\o$ and $\l$ are related by $\l= \tan\frac{\o}{2}$. 

There is an analogue of Permutability Theorem for GSGE: 

\bthm Let $\phi_i:M_0\to M_i$ be B\"acklund transformations for $n$-dimensional submanifolds in $\R^{2n-1}$ with constant $\o_i$ for $i=1, 2$.  If $\sin^2\o_1\not= \sin^2\o_2$, then there exist unique $M_3$ and B\"acklund transformations $\psi_1:M_2\to M_3$ and $\psi_2:M_1\to M_3$ with angles $\o_1, \o_2$ respectively such that $\psi_1\circ \phi_2= \psi_2\circ \phi_1$.  Moreover, if $A_i$ is the solution of the GSGE corresponding to $M_i$ for $i=0, 1, 2, 3$, then 
\beq\label{co}
A_3A_0^{-1}= (-D_2+ D_1A_2A_1^{-1}) (D_1- D_2A_2 A_1^{-1})^{-1} \I_{1, n-1},
\eeq
where $D_i= \diag(\csc\o_i, \cot\o_i, \ldots, \cot\o_i)$. 
\ethm

In other words, given a solution $A_0$ of the GSGE, we solve BT$_{A_0, \l_i}$ with $\l_i = \csc\o_i + \cot\o_i$ to get $A_i$ for $i=1, 2$. Then $A_3$ defined by the algebraic formula \eqref{co} is a solution of BT$_{A_1, \l_2}$ and $BT_{A_2, \l_1}$.  Since the constant map $A=\I$ is a solution of the GSGE, we can apply BT and permutability formula to construct infinitely many families of explicit solutions of the GSGE.

\bs
\section{Sphere congruences and  Ribaucour transforms}  \label{de}

We review the notion of sphere congruences, Christoffel and Ribaucour transforms for isothermic surfaces in $\R^3$ (cf. \cite{Da1899}). 

A {\it sphere congruence\/} in $\R^3$  is a smooth $2$-parameter family of $2$-spheres in $\R^3$:
$$S(x)=\{ c(x) +r(x) y\n y\in S^2\}, \quad x\in \co, $$
where $c:\co\to \R^3$ and $r:\co\to (0, \infty)$ are smooth maps, and $\co$ is an open subset of $\R^2$.  
A surface $f:\co\to \R^3$ is called an {\it envelope\/} of the sphere congruence $S$ if $f(p)\in S(p)$ and $f$ is tangent to the sphere $S(p)$ at $f(p)$.  To construct envelopes of $S$, we need to find a map $y:\co\to S^2$ such that $f(x)= c(x) + r(x) y(x)$ satisfying
\beq\label{at}
f_{x_1}\cdot y= f_{x_2}\cdot y=0.
\eeq
Generically there are exactly two envelopes.  
If $M$ and $\ti M$ are two envelopes of the sphere congruence $S$, then there is a natural map $\phi:M\to \ti M$ such that for each $p\in M$, there exists $x\in \co$ such that the sphere $S(x)$ is tangent to  $M$ and $\ti M$ at $p$ and $\phi(p)$ respectively. Note that the map $\phi$ determines the sphere congruence $S$.  Hence we make the following definition:

\bdefn {\bf Ribaucour transform for surfaces in $\R^3$} \hfil

A diffeomorphism $\phi:M\to \ti M$ is called a {\it sphere congruence\/} if for each $p\in M$, the normal line of $M$ at $p$ intersects the normal line of $\ti M$ at $\phi(p)$ at equal distance $r(p)$.
A sphere congruence $\phi$ from a surface $M$ in $\R^3$ to a surface $\ti M$ in $\R^3$ is called a {\it Ribaucour transform\/} if $\phi$ maps line of curvature coordinates of $M$ to those of $\ti M$. 
  \edefn

\ss\ni{\bf Isothermic surfaces}\par

An immersion $f(x_1, x_2)\in \R^3$ is called {\it isothermic\/} if  $(x_1,x_2)$ is
 both a conformal and line of curvature coordinate system. In other words, $f$ is isothermic if fundamental forms for $f$ are
\begin{equation}\label{px}
\I= e^{2q}(\rd x_1^2+ \rd x_2^2), \quad \II= e^q(r_1 \rd x_1^2 + r_2 \rd x_2^2),
\end{equation}
for some smooth functions $q, r_1$ and $r_2$.  

Set
$$w_1= e^q \rd x_1, \quad w_2= e^q \rd x_2, \quad w_{13}= r_1 \rd x_1, \quad w_{23} = r_2 \rd x_2.$$
By \eqref{ae}, $w_{12}= q_{x_2} \rd x_1 - q_{x_1} \rd x_2$.  The Gauss-Codazzi equation is:
\beq\label{am}
\bca q_{x_1x_1}+ q_{x_2x_2} + r_1r_2=0,&\\
(r_1)_{x_2}= q_{x_2} r_2, \\
(r_2)_{x_1}= q_{x_1} r_1.\eca
\eeq

For example, constant mean curvature surfaces in $\R^3$ away from umbilic points are isothermic. 

\ss\ni {\bf Ribaucour transform for isothermic surfaces}\par

Given an isothermic surface $M$ in $\R^3$, there exist an one parameter family of isothermic surfaces $M_\l$ and  Ribaucour transforms $\phi_\l : M\to M_\l$. Moreover, $\phi_\l$ can be constructed by solving a system of compatible ODEs.
  Bianchi  proved a permutability formula for these Ribaucour transforms between isothermic surfaces. 

\ss
\ni{\bf Christoffel Transform}

A {\it Christoffel transform\/}  is an orientation reversing conformal diffeomorphism $\phi:M\to \ti M$ such that $TM_p$ is parallel to $T\ti M_{\phi(p)}$ for all $p\in M$. We call $(M, \ti M)$ a {\it Christoffel pair\/}.   
Note that if $(q, r_1, r_2)$ is a solution of \eqref{am} then so is $(-q, r_1, -r_2)$.  This fact gives the Christoffel transform for isothermic surfaces:

\bthm\label{dd}
A surface $M$ in $\R^3$ is isothermic if and only if there exist a second surface $\ti M$ and a Christoffel transform $\phi:M\to \ti M$. 
Moreover, if $f(x_1, x_2)\mapsto \ti f(x_1, x_2)$ is a Christoffel transform, then the fundamental forms of $M$ and $\ti M$ are of the forms
\begin{align*}
&\I= e^{2q} (\rd x_1^2 + \rd x_2^2), \qquad \II= e^q(r_1\rd x_1^2 + r_2 \rd x_2^2),\\
&\ti \I= e^{-2q}(\rd x_1^2 + \rd x_2^2), \qquad \ti \II= e^{-q}(r_1 \rd x_1^2 - r_2 \rd x_2^2).
\end{align*}  
for some smooth solution $(q, r_1, r_2)$ of \eqref{am}. 
\ethm

\ms\ni {\bf Associated family of Christoffel pairs}

If $(f_1, f_2)$ is a Christoffel pair of isothermic surfaces in $\R^3$, then 
$$\{(\l f_1, \l f_2)\n \l\in \R\}$$ is an associated family of Christoffel pairs of isothermic surfaces in $\R^3$. The induced action of $\R^+$ on the space of solutions of \eqref{am} is
$$s\ast (q, r_1, r_2)= (q +\ln s, r_1, r_2), \quad s\in \R^+.$$

\bs
\section{Combescure transforms, O-surfaces, and $k$-tuples} 

We review the notions of conjugate coordinates on surfaces in $\R^3$, the Combescure transforms between surfaces in $\R^3$, O surfaces defined in \cite{SK03}, and $k$-tuples of k-submanifolds in $\R^n$ defined in \cite{BDPT02, DonTer08a}.

In classical geometry, a local coordinate system $(x_1, x_2)$ on a surface $M$ in $\R^3$ is said to be {\it conjugate\/} if the position function $f(x_1, x_2)$ satisfies
 $$f_{x_1x_2}= h_1 f_{x_1} + h_2 f_{x_2}$$
 for some smooth functions $h_1, h_2$; or equivalently, $\II$ is diagonalized with respect to $(x_1, x_2)$. The collection of coordinate curves $\{x_i=c_i\n c_i\in \R, i=1,2\}$ is called a {\it conjugate net} on $M$.  An orthogonal conjugate coordinate system on a surface in $\R^3$ is a line of curvature coordinate system, and the corresponding net is called an {\it O-net} (cf. \cite{Ei62}).  Note that a surface away from umbilic points admits line of curvature coordinates. 
 
  Given surfaces $M, \ti M$ in $\R^3$, 
  a diffeomorphism $\phi:M\to \ti M$ is a {\it Combescure transform\/} if $TM_p= T\ti M_{\phi(p)}$ for all $p\in M$.
  These classical notions can be generalized to submanifolds in Euclidean spaces as follows:
  
  \ss\ni {\bf Conjugate coordinate system for submanifolds in $\R^n$}
  
   A coordinate system $x$ on a $k$-dimensional submanifold $M$ in $\R^n$ is called {\it conjugate\/} if the position function $f(x)$ satisfies the following conditions:
  $$f_{x_i x_j} = \sum_{\ell=1}^k c_{ij\ell} f_{x_\ell}, \quad 1\leq i < j \leq k$$
  for some smooth functions $c_{ij\ell}$.  We call the collection of all coordinate curves of a conjugate coordinate system a {\it conjugate net\/} on the submanifold.  
  
  If $f(x)$ is an immersion parametrized by conjugate coordinate system, then  $f_{x_i}$ are eigenvectors of the shape operator $A_v$ along any normal vector field $v$. So all shape operators commute, which implies that the normal bundle of $f$ must be flat.  An orthogonal conjugate coordinate system on a submanifold in $\R^n$ is a line of curvature coordinate system. Unlike surfaces in $\R^3$, submanifolds in Euclidean space with flat normal bundle generically do not admit line of curvature coordinate systems.
  
  \bdefn {\bf Combescure transform for submanifolds} \hfil
  
   A diffeomorphism $\phi$ from a $k$-dimensional submanifold $M$ to another $\ti M$ in $\R^n$ is called a {\it Combescure transform\/} if $TM_p= T\ti M_{\phi(p)}$ for all $p\in M$.   
   \edefn
   
   \bdefn {\bf Combescure O-transform} \cite{DonTer08a}\hfil
   
 Let $M, \ti M$ be submanifolds in $\R^n$ admitting line of curvature coordinates (so they have flat normal bundles).  A Combescure transform $\phi:M\to \ti M$ is called a {\it Combescure O-transform\/} if 
\ben
\item $\phi$ preserves line of curvature coordinates,
\item if $v$ is parallel normal field on $M$, then $v$ is a parallel normal field on $\ti M$ (since $TM_p= T\ti M_{\phi(p)}$ for all $p\in M$, we can identify $\nu(M)_p$ as $\nu(\ti M)_{\phi(p)}$).
\een
 \edefn
  
\bdefn \label{ei} {\bf Combescure O-map} \cite{DonTer08a} \hfil

Let $\W$ be an open subset of $\R^k$, and $\calM_{n\times \ell}$ the space of real $n\times \ell$ matrices with $\ell\leq k$.  A smooth map $Y= (Y_1, \ldots, Y_\ell):\W\to \calM_{n\times \ell}$ is called a {\it Combescure O-map\/} if it satisfies the following conditions:
 \ben
 \item[(a)] Each $Y_i:\W\to \R^n$ is an immersion with flat normal bundle and parametrized by line of curvature coordinates.
 \item[(b)] The map $Y_i(x)\mapsto Y_{i+1}(x)$ is a Combescure O-transform for $1\leq i\leq \ell-1$.
\item[(c)] Let  $e_i$ be the unit direction of $(Y_1)_{x_j}$ for $1\leq j\leq k$ (so $e_i$ is parallel to $(Y_i)_{x_i}$ for $2\leq i\leq \ell$), and $a_{ij}$'s defined by 
  $(Y_i)_{x_j} = a_{ij}e_j$ for  $1\leq i\leq \ell$ and $j\leq k$.   We call $(a_{ij})$ the {\it metric matrix associated to $Y$}. The rank of $(a_{ij}(x))$ is $\ell$ for all $x\in \W$, 
 \een
\edefn

\brem\label{ed}
  Let $Y=(Y_1, \ldots, Y_\ell):\W\to \calM_{n\times \ell}$ be a Combescure O-map, and $(a_{ij})$ the metric matrix associated to $Y$. 
Let $(e_{k+1}, \ldots, e_n)$ be an orthonormal parallel normal frame for $Y_1$, and $g=(e_1, \ldots, e_n)$. Then $g$ is an adapted frame on $Y_j$ for all $1\leq j\leq k$.  Hence they have the same Maurer-Cartan form $g^{-1}\rd g= (w_{AB})$. 
 By Cartan Lemma \ref{fi} and \eqref{ae}, we have
  $$w_{rs}= \frac{(a_{ir})_{x_s}}{a_{is}} \rd x_r - \frac{(a_{is})_{x_r}}{a_{ir}} \rd x_s, \quad 1\leq r\not= s\leq k, \quad1\leq i\leq \ell.$$
  So $$\frac{(a_{ir})_{x_s}}{a_{is}}= \frac{(a_{1r})_{x_s}}{a_{1s}}, \quad 1\leq r\not= s \leq k, \, 1\leq i\leq \ell.$$ 
  Geometrically, this means that $\K_j e_i=\K_1 e_i$ for all $i, j\leq k$, where $\K_j$ is the Levi-Civita connection of the induced metric $\I_j$ of $Y_j$.  
  Since $x$ is a line of curvature coordinate system, there exist smooth functions $h_{i\a}$ such that
  $$w_{i\a}= h_{i\a} \rd x_i, \quad 1\leq i\leq k, k< \a\leq n.$$    \erem
   
   \bdefn  {\bf O surfaces} (\cite{SK03})
   
   Two surfaces $f_1(x), f_2(x)$ in $\R^3$ parametrized by line of curvature coordinates are called {\it O-surfaces\/} if 
   \ben
   \item[(a)] the map $f_1(x)\mapsto f_2(x)$ is a Combescure transform for all $i\not=j$,
   \item[(b)] $\frac{(a_{21})_{x_2}}{a_{22}} = \frac{(a_{11})_{x_2}}{a_{12}}$ and $ \frac{(a_{22})_{x_1}}{a_{21}}= \frac{(a_{12})_{x_1}}{a_{11}}$,  where $e_j$ is the unit direction of $(f_1)_{x_j}$ (hence $e_j$  is parallel to $(f_2)_{x_j}$) for $j=1,2$ and $a_{ij}$'s  are defined by 
   $(f_i)_{x_j}= a_{ij} e_j$ for  $i,j=1,2$.
\een
\edefn

As a consequence of Remark \ref{ed}, we have

\bprop Two surfaces $f(x), \ti f(x)$ parametrized by line of curvature coordinates are O surfaces if and only if the map $(f, \ti f)$ is a Combescure O-map. \eprop

\bdefn {\bf $k$-tuples in $\R^n$} \cite{BDPT02, DonTer08a} \hfil\par

A Combescure O-map $Y=(Y_1, \ldots, Y_k)$ of $k$-dimensional submanifolds in $\R^n$ is called 
\ben
\item a {\it $k$-tuple of $k$-submanifolds in $\R^n$ of type $\R^{k-\ell, \ell}$} (or just $k$-tuple in $\R^n$ of type $\R^{k-\ell, \ell}$) if all rows of the metric matrix of $Y$ have constant length in $\R^{k-\ell, \ell}$. 
\item a {\it $k$-tuples of $k$-submanifolds in $\R^n$ of type $O(k-\ell, \ell)$} if the metric matrix of $Y$ lies in $O(k-\ell, \ell)$,
\item a a {\it $k$-tuple of $k$-submanifolds in $\R^n$ of null $\R^{k-\ell, \ell}$ type} if all rows of the metric matrix of $Y$ are null vectors in $\R^{k-\ell, \ell}$.
\een
\edefn

Combescure O-maps and $2$-tuples occur naturally in surface geometry:

\beg\label{ej1} 
If $f(x)$ is a surface in $\R^3$ parametrized by line of curvature coordinates, then $\{f, e_3\}$ and $\{f, f+ re_3\}$ are O surfaces in $\R^3$, where $r\in \R$ is a constant and $e_3$ is the unit normal.   
\eeg

\beg\label{ej2}
 A Christoffel pair of isothermic surfaces $(f_1, f_2)$ is a Combescure O-map whose metric matrix $(a_{ij})$ is of the form $\bpm e^q& e^q\\ e^{-q} & -e^{-q}\epm$ for some $q$, i.e., it is a $2$-tuple in $\R^3$ of null $\R^{1, 1}$ type.
 \eeg
 
 \beg \label{ej3} \cite{BDPT02, SK03}
 
 A $2$-tuple $(f_1, f_2)$ of surfaces in $\R^3$ of type $O(1,1)$ is a Combescure O-map whose metric matrix is of the form $\bpm \cosh q&\sinh q\\ \sinh q & \cosh q\epm$, and the two fundamental forms for $Y_1, Y_2$ are  
\begin{align*}
&\bca \I_1= \cosh^2 u \rd x_1^2 + \sinh^2 u \rd x_2^2, \\ \II_1=r_1\cosh u \rd x_1^2 + r_2 \sinh u \rd x_2^2,\eca, \\
& \bca \I_2= \sinh^2 u \rd x_1^2 + \cosh^2 u \rd x_2^2, \\ \II_2= r_1 \sinh u\rd x_1^2 + r_2 \cosh u \rd x_2^2.\eca 
\end{align*}
Note that 
\ben
\item the Gaussian curvature of $f_1$ and $f_2$ are equal, $K_1(x)= K_2(x)$,
\item  $(Y_1+Y_2, Y_1-Y_2)$ is an isothermic pair.
\een
\eeg

\beg\label{ej4} \cite{BDPT02, SK03}
 A $2$-tuple $(f_1, f_2)$ of surfaces in $\R^3$ of type $O(2)$ is a Combescure O-map whose metric matrix is of the form $\bpm \cos q & \sin q\\ -\sin q & \cos q\epm$, and the fundamental forms of $Y_1, Y_2$ are
$$\bca \I_1= \cos^2 q\rd x_1^2  +\sin^2 q\rd x_2^2,\\ \II_1= r_1 \cos q \rd x_1^2 + r_2 \sin q \rd x_2^2,\eca \quad \bca \I_2= \sin^2 q \rd x_1^2 + \cos^2 q \rd x_2^2, \\ \II_2= r_1 \sin q\rd x_1^2 - r_2 \cos q \rd x_2^2.
\eca$$
Thus the Gaussian curvature $K_1(x)= -K_2(x)$.  

If $f_1(x)$ is a surface with $K=-1$ parametrized by Tchebyshef line of curvature coordinates as in section \ref{cz}, then $r_1= \sin q$, $r_2= -\cos q$, $q$ is a solution of SGE, and $(f, e_3)$ is a $2$-tuple of surfaces in $\R^3$ of type $O(2)$.
\eeg

\bdefn \label{fg} {\bf Isothermic$_\ell$ $k$-submanifolds in $\R^n$} \cite{DonTer08a}\hfil\par

A $k$-dimensional submanifold $M$ in $\R^n$ is {\it isothermic$_\ell$} if 
\ben
\item the normal bundle is flat,
\item there is a line of curvature coordinate system $(x_1, \ldots, x_k)$ such that $\I=\sum_{i=1}^k g_{ii} \rd x_i^2$ satisfies $\sum_{i=1}^{k-\ell} g_{ii} -\sum_{i= k-\ell +1}^kg_{ii}=0$.  
\een
\edefn

\brem\hfil\par

\ben
\item A $k$-tuple in $\R^n$ of type $O(k-\ell, \ell)$ is of type $\R^{k-\ell, \ell}$,
\item a Christoffel pair of isothermic surfaces in $\R^n$ is a $2$-tuple in $\R^n$ of null $\R^{1,1}$ type (cf. \cite{BDPT02, Bu04})    
\item The equation for $k$-tuples in $\R^n$ of type $\R^{k-\ell, \ell}$ is the $\frac{O(n+k-\ell, \ell)}{O(n)\times O(k-\ell, \ell)}$-system and there are Ribaucour transforms and permutability formulae for these $k$-tuples. These will be reviewed in sections 8 and 10.
\item If $Y=(Y_1, \ldots, Y_k)$ is a $k$-tuple in $\R^n$ of null $\R^{k-\ell, \ell}$ type, then each $Y_i$ is an isothermic$_\ell$ submanifold in $\R^n$ and $Y_i$ and $Y_j$ are related by Combescure $O$-transforms.  
\een  
\erem

\bs
\section{From moving frame to Lax pair} \label{da}

Suppose the PDE for $q:\R^n\to V$ has a $\calG$-valued {\it Lax pair\/} $\o_\l$ on $\R^n$, where $\calG$ is the Lie algebra of a Lie group $G$.  If $q$ is a solution of the PDE, then given $c_0\in G$ there is a unique $G$-valued solution $E(x,\l)$ for 
$$E^{-1} \rd E= \o_\l, \quad E(0,\l)=c_0,$$
which will be called a {\it parallel frame} of the solution $q$ or of its Lax pair $\o_\l$. The solution with initial data $c_0=\I$ is called the {\it normalized parallel frame\/}.  

The existence of a Lax pair is one of the characteristic properties of soliton equations. The SGE, GSGE, and the Gauss-Codazzi equation for isothermic surfaces and for flat Lagrangian submanifolds in $\C^n$, and the equation for $k$-tuples in $\R^n$ of type $\R^{k-\ell, \ell}$ are soliton equations and their Lax pairs were found in \cite{AKNS74, ABT86, CGS95, TW08, DonTer08a} respectively.  In general, it is not easy to determine whether a PDE has a Lax pair.  We explain in this section how to construct
\ben
\item Lax pairs for SGE, GSGE, equations for flat Lagrangian submanifolds in $\C^n$, and for $k$-tuples in $\R^n$ of type $\R^{k-\ell, \ell}$  from the Maurer-Cartan forms of specially chosen moving frames of the {\it associated family\/} of these submanifolds,
\item  the immersions of these submanifolds from parallel frames of the corresponding Lax pairs.
\een

\ss\ni{\bf $K=-1$ surfaces in $\R^3$}

\ss\ni {\it Lax pair} (cf. \cite{Bo94, TU00b})

Suppose $M$ is a surface in $\R^3$ with $K=-1$, $(s,t)$ the Tchebyshef asymptotic coordinate system, and $q(s,t)$ is the solution of SGE corresponding to $M$.  Let $f^\l:M^\l\to \R^3$ denote the $K=-1$ surface corresponding to the solution  $q^\l(s,t)= q(\l s, \l^{-1}t) $.  We derive a Lax pair for the SGE from the Maurer-Cartan form for $M^\l$: For each non-zero real $\l$, choose the orthonormal frame $F^\l= (e_1^\l, e_2^\l, e_3^\l)$  on $M^\l$ such that 
$e_1^\l=f^\l_{s}$, and $e_3^\l$ is the unit normal to $M^\l$. Set $w^\l:= (F^\l)^{-1} \rd F^\l$. Substitute $(\frac{1}{2\l} s, 2\l t)$ for $(s,t)$ in $\w^\l$ to get a one-parameter family of flat $o(3)$-valued connection $1$-forms:
\beq\label{ci}
\w^\l= \bpm 0& -2q_s&0\\ 2q_s&0 & -2\l\\ 0 & 2\l &0\epm \rd s + \frac{1}{2\l}\bpm 0&0&\sin 2q\\ 0&0&\cos 2q\\ -\sin 2q & -\cos 2q & 0\epm \rd t.
\eeq
  To get the known Lax pair of SGE, we identify the Lie algebra $o(3)$ as $su(2)$ to rewrite the family $\w^\l$ of $o(3)$-valued connections  as a family of flat $su(2)$-valued connection $1$-forms:
\beq\label{bb}
\o_\l= \bpm -\im\l & -q_s\\ q_s & \im\l\epm \rd s + \frac{\im}{4\l} \bpm \cos 2q& -\sin 2q\\ -\sin 2q & -\cos 2q \epm \rd t.
\eeq
Moreover, given $q:\R^2\to \R$, then $q$ is a solution of the SGE if and only if $\o_\l$ defined by \eqref{bb} is flat for all non-zero $\l\in \C$.

\ss\ni {\it Sym's formula} \cite{Sym85}

If $q$ is a solution of the SGE, then we can construct the corresponding surface with $K=-1$ in $\R^3$ from a parallel frame of the Lax pair associated to $q$ as follows: 
Set $\o_\l$ by \eqref{bb}, and let $E(s, t,\l)$ be a parallel frame for $\o_\l$, i.e., the solution of 
$$E^{-1}\rd E= \o_\l, \quad E(0, 0, \l)=c_0\in SU(2).$$
Since $\o_{\bar\l}^*+\o_\l=0$, $E(s,t,\bar\l)^*E(s,t,\l)= \I$.  Hence $E(s,t, r)\in SU(2)$ for any real number $r$. 
Set $$f_r= \frac{\p E}{\p \l} E^{-1}\,\bigg|_{\l=r}.$$ Because $E(s, t, r)\in SU(2)$ for $r\in \R$, we have $f_r\in su(2)$.  Also
$$\rd f_r = E(s,t,r)\left( \bpm -\im &0\\ 0& \im\epm \rd s +\frac{\im}{4r^2}\bpm -\cos 2q & \sin 2q\\ \sin 2q & \cos 2q\epm\rd t\right) E(s,t,r)^{-1}.$$ 
If we identify $su(2)$ as $\R^3$, then
$f:=f_{\frac{1}{2}}(s,t)$ is a surface with $K=-1$, $(s,t)$ is the Tchebyshef asymptotic coordinate system,\ and $q$ is the solution of the SGE corresponding to $f$.

\ss\ni{\bf $n$-submanifolds in $\R^{2n-1}$ with sectional curvature $-1$}

\ss\ni {\it Lax pair}

Let $f:M^n\to \R^{2n-1}$ be an immersion with sectional curvature $-1$, $x$ the Tchebyshef line of curvature coordinate system, $e_i$ the unit direction of $f_{x_i}$, $(e_{n+1}, \ldots, e_{2n-1})$ the parallel normal frame,  and 
$$\I=\sum_{i=1}^n a_{1i} \rd x_i^2, \quad \II= \sum_{i=1, j=2}^n a_{1i} a_{ji} \rd x_i^2 e_{n+j-1},$$
the fundamental forms as in Theorem \ref{bn}. Set
$F=(f_{ij})$ as in \eqref{bo2}, $w_i=a_{1i} \rd x_i$, $ w_{i, n+j-1}= a_{ji} \rd x_i$, $w_{ij}= f_{ij} \rd x_i - f_{ji} \rd x_j$ and $ w_{n+i-1, n+j-1}=0$.

We associate to the immersion $f$ two flat connection $1$-forms:
The sectional curvature of $\I=\sum_{i=1}^n w_i^2$ 
is $-1$, giving  $$\rd w_{ij} + \sum_k w_{ik}w_{kj} = - w_i w_j,$$ which is equivalent to 
$$\zeta_1= \bpm \w & \xi^t\\ \xi &0\epm, \qquad \w = (w_{ij})_{i, j\leq n}, \quad \xi= (w_1, \ldots, w_n)$$
being a flat $o(n,1)$-valued connection $1$-form.  The Maurer-Cartan form of $f$ gives a flat $o(2n-1)$-valued $1$-form 
$$\varpi= g^{-1}\rd g= (w_{AB})= \bpm w & \eta\\ -\eta^t & 0\epm,$$ where $g=(e_1, \ldots, e_{2n-1})$,
$ w=(w_{ij})_{i,j\leq n}$ and $ \eta_{ij}= a_{ji} \rd x_i$.   

It is easy to see that an $o(2n-1)$-valued $1$-form $\bpm \w & \eta\\ -\eta^t &0\epm$ is flat if and only if 
$$\zeta_2=\bpm \w & \im\eta \\ \im\eta^t&0\epm= \bpm 1&0\\ 0& -\im\epm \bpm w&\eta\\ -\eta^t & 0\epm \bpm 1 & 0\\ 0& \im\epm$$ is a flat $o(n, n-1, \C)$-valued $1$-form.  We embed $o(n,1)$ and $o(n, n-1)$ into $o(n,n)$ as Lie subalgebras by
\begin{align*}
&o(n,1)=\{y=(y_{ij})\in o(n,n)\n y_{ij}=0, \,\, \forall\,\, n+1< i, j\leq 2n\}, \\
& o(n,n-1)= \{y=(y_{ij})\in o(n,n)\n y_{i, n+1}= y_{n+1, i}=0 \, \, \forall \,\, 1\leq i\leq n\}.
\end{align*} 
Use these embeddings to write $\zeta_1, \zeta_2$ as flat $o(n,n)$-valued $1$-forms: 
$$\zeta_j= \bpm \w & \d A^t D_j\\ D_j A\d &0\epm,\qquad j=1, 2,$$
where $\d=\diag(\rd x_1, \ldots, \rd x_n)$,  $D_1= e_{11}=\diag(1, 0,\ldots, 0)$ and $D_2= \im(\I-e_{11})$.   
The flatness of $\zeta_1$ and $\zeta_2$ gives:
\beq\label{bc}
\bca \rd \w+ \w\wedge \w+ \d A^t D^2 A\d =0, & \quad {\rm where\,} w= \d F- F^t\d, \\
D\rd A \wedge \d + DA\d \wedge \w=0\eca
\eeq
for $D= D_1$ or $D= D_2$.  Write system \eqref{bc} in terms of $A$ and $F$ to get the GSGE \eqref{bp}.

Set $D_\o=\cos\o D_1 + \sin \o  D_2=\diag(\cos\o, \im\sin\o, \ldots, \im\sin\o)$. Then $D_\o^2= -\sin^2\o \I + D_1$.  Since $\d \wedge \d=0$,  $\d A^tD_\o^2 A\d = \d A^t e_{11} A\d$.  So \eqref{bc} is flat for all $D= \frac{1}{2}(e^{i\o} \I + e^{-i\o} \I_{1, n-1})$. 
Hence 
\beq\label{bs}
\o_\l = \bpm \w& \d A^t D_\l\\ D_\l A\d& 0\epm,
\eeq
is a flat $o(n,n)$-valued connection $1$-form on $\R^n$ for all $\l= e^{i\o}$, where $D_\l = \frac{1}{2} (\l \I+ \l^{-1}J)$, $ w=\d F- F^t \d$, and $ A^{-1}\rd A= \d F^t- F\d$.  Moreover, $A$ is a solution of GSGE if and only if $\o_\l$ is flat for all $\l\not=0$.  This is the Lax pair given in \cite{ABT86} for the GSGE.

\ss\ni {\it SGE has two Lax pairs}

 Note that SGE has two Lax pairs, one is the $sl(2,\C)$-valued connection $1$-form \eqref{bb} in asymptotic coordinates and the other is the $o(2,2)$-valued connection $1$-form \eqref{bs} in line of curvature coordinates.

\ss\ni {\it Construct immersions}
 
Suppose $(A=(a_{ij}), F=(f_{ij}))$ is a solution of the GSGE \eqref{bp2}, and $\o_\l$ the Lax pair defined by \eqref{bs}.  Let $E(x,\l)$ denote the normalized parallel frame of $\o_\l$, and
$$g(x):=\bpm 1& 0\\ 0& \im\epm E(x,\im) \bpm 1& 0\\ 0& -\im\epm.$$
Then $g(x)\in O(2n)$, $g(x)_{n+1, i}= g(x)_{i, n+1}=0$ for $i\not= n+1$, $g(x)_{n+1, n+1}=1$, and
$$g^{-1}\rd g= \bpm 1 &0\\ 0 & \im\epm \o_\im \bpm 1 &0\\ 0& -\im\epm = \bpm w & \d A^t(\I-e_{11})\\ -(\I-e_{11}) A\d &0\epm$$
 is a flat $o(2n-1)$-valued connection $1$-form with $g(0)=\I$. 
  Hence $g(x)\in O(2n-1)$. Let $e_i(x)$ denote the $i$-th column of $g(x)$. Then the following system  
\beq\label{cl}
\rd f= \sum_{i=1}^n a_{1i}e_i\rd x_i
\eeq
is solvable for $f$ in $\R^{2n-1}$ and the solution $f$ (up to translation) has sectional curvature $-1$.

\ss\ni{\bf Flat Lagrangian submanifolds in $\C^n$} \cite{TW08}

\ss\ni {\it Egoroff line of curvature coordinate system}

If $f:M\to \C^n=\R^{2n}$ is a flat Lagrangian submanifold with flat and non-degenerate normal bundle, then there exist a coordinate system $(x_1, \ldots, x_n)$ and function $\phi$ such that
\beq\label{cn}
\bca\I= \sum_{i=1}^n\phi_{x_i} \rd x_i^2,\\ \II=\sum_{i=1}^n \rd x_i^2 \otimes J(f_{x_i}),\eca
\eeq
where $J$ is the standard complex structure on $\R^{2n}$. We call $x$ the {\it Egoroff line of curvature coordinate system\/}.  
Let $$e_i= \frac{f_{x_i}}{(\phi_{x_i})^{\frac{1}{2}}}, \quad e_{n+i} = Je_i, \quad  1\leq i\leq n,$$ $g=(e_1, \ldots, e_{2n})$ the adapted frame for $f$, and $\varpi= g^{-1}\rd g = (w_{AB})$.
Then the dual $1$-forms for $e_1, \ldots, e_n$ are $w_i= (\phi_{x_i})^{\frac{1}{2}} \rd x_i$ and $w_{i, n+j} =\d_{ij} \rd x_i$.  By the Cartan Lemma \ref{fi} and \eqref{ae}, we have 
$$w_{ij}= \b_{ij} \rd x_i - \b_{ji} \rd x_j, \quad {\rm where\,\,}  \b_{ij}=\bca \frac{\phi_{x_ix_j}}{2(\phi_{x_i}\phi_{x_j})^{\frac{1}{2}}}, & i\not=j,\\ 0, &i= j\eca$$
for $i, j\leq n$.  
Note that $\b=(\b_{ij})$ is symmetric.  Set
 $h=(\phi_{x_1}^{\frac{1}{2}}, \ldots, \phi_{x_n}^{\frac{1}{2}})^t$.
The Gauss-Codazzi equation and the structure equation for $f$ is the PDE for $(\b, h)$ defined by the condition that 
\beq\label{dv}
\tau=\bpm [\d, \b]&  \d & \d h\\ - \d & [\d,\b] &0\\ 0&0&0\epm
\eeq
is flat, i.e., it is the following system for $(\b, h)$:
\beq \label{dh}
\bca (h_i)_{x_j}= \b_{ij} h_j, & i\not=j,\\
(\b_{ij})_{x_i} + (\b_{ij})_{x_j} + \sum_k \b_{ik} \b_{kj}=0, & i\not=j,\\
(\b_{ij})_{x_k}= \b_{ik} \b_{kj}, & i, j, k \quad {\rm distinct}.\eca
\eeq
Conversely, if $(\b, h)$ is a solution of \eqref{dh} with $\b=(\b_{ij})$ symmetric and $h=(h_1, \ldots, h_n)^t$, then the first equation of \eqref{dh} implies that $h_i(h_i)_{x_j}=h_j(h_j)_{x_i}$ for all $i\not= j$.  So $(h_1^2, \ldots, h_n^2)$ is a gradient field, i.e., there is a function $\phi$ such that $h_i^2= \phi_{x_i}$ for $1\leq i\leq n$.  Hence there is a flat Lagrangian immersion $f(x)$ in $\C^n$ such that $\I, \II$ are of the form \eqref{cn}. 

\ss\ni {\it Associated family of flat Lagrangian submanifolds in $\C^n$}

If $M$ is a flat Lagrangian submanifold in $\C^n$ with $\I, \II$ as in \eqref{cn}, then given $\l\in \R$, there is a flat Lagrangian submanifold $M_\l$ in $C^n$ with
$$\I_\l = \I = \sum_{i=1}^n\phi_{x_i} \rd x_i^2, \qquad \II_\l = \l\sum_{i=1}^n \rd x_i^2 \otimes J(f_{x_i}).$$

\ss\ni {\it Lax pair}

 If $f_\l$ is the associated family of $f$, then the Maurer-Cartan form \eqref{dv} for $f_\l$ is 
\beq \label{dm}
 \o_\l = \bpm [\d, \b] & \l \d & \d h\\ -\l \d &[\d, \b] & 0\\ 0 & 0 & 0\epm.
\eeq
Moreover, the following statements are equivalent: (i) $\o_1$ is flat, (ii) $\o_\l$ is flat for all $\l\in \C$, (iii) $(\b, h)$ is a solution of \eqref{dh}.  

\ss\ni {\it Construct flat Lagrangian submanifold from parallel frame}

If $(\b, h)$ is a solution of the \eqref{dh} and $E$ a parallel frame of $\o_\l$ given by \eqref{dm}, then
\ben
\item there exists $\phi$ such that $h_i^2= \phi_{x_i}$ for $1\leq i\leq n$,
\item for each real $r$, $E(x,r)$ is of the form $\bpm g(x,r)& f(x,r)\\ 0 &1\epm$ with $g\in U(n)\subset O(2n)$ and $f(\cdot, r)\in \R^{2n}$, 
\item $f(\cdot, r)$ is a flat Lagrangian submanifold in $\C^n=\R^{2n}$ with 
$$\I_r=\sum_{i=1}^n \phi_{x_i} e_i\rd x_i^2, \quad \II_r= \sum_{i=1}^n r \rd x_i ^2 \otimes J(f_{x_i}(\cdot, r)),$$
where $e_i(x,r)$ is the $i$-th column of $g(x,r)\in O(2n)$ and $J(f_{x_i})$ is parallel to $e_{n+i}$ for $1\leq i\leq n$
\een

\ss\ni {\bf Isothermic surfaces in $\R^3$} 

\ss\ni {\it Lax pair}

We use the associated family of  Christoffel transforms to construct a Lax pair for isothermic surfaces. 
Suppose $(f(x), \ti f(x))$ is a Christoffel transform of isothermic immersions in $\R^3$. Let $e_1, e_2$ denote the coordinate directions. By Theorem \ref{dd}, there exists a solution $(q, r_1, r_2)$ of \eqref{am} such that 
$$\rd f= e^q(\rd x_1 e_1 + \rd x_2 e_2), \qquad \rd \ti f= e^{-q} (\rd x_1 e_1 - \rd x_2 e_2).$$
Write the above equation in matrix form:
$$\rd (f, \ti f)= (e_1, e_2, e_3)\bpm \rd x_1 & 0\\ 0 & \rd x_2 \\ 0&0\epm \bpm \cosh q& \sinh q\\ \sinh q & \cosh q\epm\bpm 1&1\\ 1& -1\epm.$$
Set 
\begin{align*}
&\zeta:= (e_1, e_2, e_3)\bpm \rd x_1 & 0\\ 0 & \rd x_2 \\ 0&0\epm \bpm \cosh q& \sinh q\\ \sinh q & \cosh q\epm, \\
& \xi=\bpm 0 & \zeta\\ -J\zeta^t &0\epm, \quad J=\diag(1, -1). 
\end{align*}
Compute directly to see that 
$\rd \xi=0$ and $\xi\wedge\xi=0$, which implies that $\xi$ is a flat $0(4,1)$-valued connection $1$-form.  
Apply the above computation to  the associated family $\l(f, \ti f)$ to see that $\l\xi$ is a flat connection $1$-form for all $\l\in \R$. 
Set 
$$g_1=(e_1, e_2, e_3), \quad g_2=\bpm \cosh q& -\sinh q\\ -\sinh q &\cosh q\epm, \quad g=\bpm g_1 &0\\ 0& g_2\epm.$$
The gauge transformation of $\l\xi$ by $g^{-1}$ is
  \beq\label{df}\o_\l = \l g^{-1}\xi g + g^{-1}\rd g = \bpm w& \l D\\ -\l J D^t & \tau\epm,
 \eeq
  where 
  \begin{align*}
  &w = g_1^{-1}\rd g_1=\bpm 0 & q_{x_2}\rd x_1-q_{x_1} \rd x_2 & r_1 \rd x_1\\ -q_{x_2}\rd x_1+q_{x_1} \rd x_2 &0& r_2 \rd x_2\\ -r_1 \rd x_1 & -r_2 \rd x_2 &0\epm, \\
  & \tau= g_2^{-1}\rd g_2= \bpm 0& -\rd q\\ -\rd q &0\epm, \quad D=\bpm \d\\ 0\epm, \quad \d=\diag(\rd x_1, \rd x_2).
  \end{align*} 
  Since $\xi\l$ is flat, so is $\o_\l$.
  Moreover, $(q, r_1, r_2)$ is a solution of the Gauss-Codazzi equation \eqref{am} of isothermic surfaces if and only if 
  $\o_\l$ is flat for all parameters $\l$.  In other words, $\o_\l$ is a Lax pair of the isothermic equation \eqref{am}.
  
  Note that $\o_\l$ can be written as 
  \beq\label{ck} 
  \o_\l=\sum_{i=1}^2 (a_i\l + [a_i, v]) \rd x_i, 
  \eeq
  where $J=\I_{1,1}$, 
 \begin{gather}
  a_i =  \bpm 0 & D_i\\ -JD_i^t & 0\epm, \quad D_1= \bpm 1&0\\ 0&0\\0&0\epm, \quad 
  D_2=\bpm 0&0 \\ 0&1\\ 0&0\epm,  \label{cj1}\\
   v=\bpm 0& \eta\\ -J \eta^t& 0\epm, \quad \eta= \bpm 0 & q_{x_1} \\ -q_{x_2} & 0\\-r_1 & r_2\epm.\label{cj2}
  \end{gather} 
  
  \ms\ni {\it Construction of Christoffel pairs of isothermic surfaces from parallel frames}
  
  \ss\ni {\it Method 1}
  
  Let $(q, r_1, r_2)$ be a solution of \eqref{am}, $\o_\l$ its Lax pair defined by \eqref{df}, and $E(x,\l)$ a parallel frame for $\o_\l$ with initial data $c_0\in O(3)\times O(1,1)$.  Since $\o_0\in o(3)\times o(1,1)$, 
  $g:= E(x,0)=\bpm g_1 &0\\ 0& g_2\epm\in O(3)\times O(1,1)$. Write 
  $$g_1=(e_1, e_2, e_3), \quad g_2=\bpm \cosh q & -\sinh q\\ -\sinh q & \cosh q\epm.$$ Then
  $$g\ast\o_\l= g\o_\l g^{-1}-\rd g g^{-1}= \l \bpm 0& \zeta\\ -J\zeta^t & 0\epm,$$
 where 
 $$\zeta=(e_1, e_2, e_3)\bpm \rd x_1 &0\\ 0 & \rd x_2\\ 0&0\epm \bpm \cosh q & \sinh q\\ \sinh q& \cosh q\epm.$$ 
 The flatness $g\ast\o_\l$ implies that $\rd \zeta=0$. Hence there exists a $3\times 2$ matrix valued map $Y$ such that 
  $\rd Y= \zeta$.  Moreover, $(f_1, f_2)= Y\bpm 1& 1\\ 1 &-1\epm$ is a Christoffel pair of isothermic surfaces in $\R^3$ and $(q, r_1, r_2)$ is the corresponding solution of \eqref{am}.  
  
  \ss\ni {\it Method 2}
  
  We claim that if $E$ is the normalized parallel frame of the Lax pair $\o_\l$ defined by \eqref{df} of a solution $(q, r_1, r_2)$ of \eqref{am}, then $\frac{\p E}{\p \l}E^{-1}\,\big|_{\l=0}$ is of the form $\bpm 0 & Z\\ -\I_{1,1} Z^t& 0\epm$ for some $3\times 2$ matrix value map $Z$ and $(f_1, f_2)= Z\bpm 1&1\\1&-1\epm$ is a Christoffel pair of isothermic surfaces in $\R^3$ with fundamental forms as in Theorem \ref{dd}.   
 
To see this, we first note that $\o_\l$ is $o(4,1, \C)$-valued $1$-form and satisfies the $\frac{O(4,1)}{O(3)\times O(1,1)}$ reality condition:
 $$\overline{\o_{\bar\l}}= \o_\l, \quad  \I_{3,2} \o_\l \I_{3,2}= \o_{-\l}.$$  So the normalized parallel frame 
  $E$ of $\o_\l$ satisfies $E(x,\l)\in O(4,1, \C)$ satisfying
  \beq\label{fe}
  \overline{E(x, \bar\l)}= E(x,\l), \quad \I_{3,2}E(x, \l) \I_{3,2}= E(x, -\l).
  \eeq 
  Note that $\tau(y)= \bar y$ and $\sigma(y)= \I_{3,2}y\I_{3,2}^{-1}$ are involutions on $O(4,1, \C)$ that give the symmetric space $\frac{O(4,1)}{O(3)\times O(1,1)}$, and 
  $$o(4,1)= \calK + \calP, \quad \calK= o(3)\times o(1,1), \, \calP= \left\{ \bpm 0 & \xi \\ -\I_{1,1}\xi^t & 0 \epm\right\}$$
is  the Cartan decomposition of $\pm 1$ eigenspaces of $\sigma$ on the fixed point set $o(4,1)$ of $\tau$.  
  It follows from \eqref{fe} that
   $\frac{\p E}{\p \l} E^{-1}\big|_{\l=0}$ lies in $\calP$, hence is of the form $\bpm 0 & Z\\ -JZ^t &0\epm$ for some $3\times 2$ valued map $Z$.  A direct computation implies that $$\rd Z= g_1 \bpm \d \\ 0\epm g_2^{-1}$$ and $Z\bpm 1& 1\\ 1& -1\epm$ is a Christoffel pair associated to the solution $(q, r_1, r_2)$, where $g_1, g_2$ is given by $E(x,0)=\bpm g_1 &0\\ 0 &  g_2\epm$.

\ms
\ni {\bf $k$-tuples in $\R^n$ of type $\R^{k-\ell, \ell}$\/} \cite{BDPT02, DonTer08a}

\ss\ni {\it Lax pair}

First we associate to a $k$-tuple $Y$ in $\R^n$ of type $\R^{k-\ell, \ell}$ two flat connections, and then use them to construct a Lax pair for the equation of $Y$. 

\bthm  \label{eg} \cite{DonTer08a}
Let $Y=(Y_1, \ldots, Y_k):\W\to \calM_{n\times k}$ be a $k$-tuple in $\R^n$ of type $\R^{k-\ell, \ell}$, $e_j$ the unit direction of $(Y_1)_{x_j}$ for $1\leq j\leq k$, $e_{k+1}, \ldots, e_n$ a parallel orthonormal normal frame for $Y_1$, $g=(e_1, \ldots, e_n)$, $(w_{ij})_{i,j\leq n}= g^{-1}\rd g$, and $(a_{ij})_{i, j\leq k}$ the metric matrix associated to $Y$ defined by
$(Y_i)_{x_j} =  a_{ij}e_j$ for $ 1\leq i, j\leq k$.  
Set
$$f_{ij} =\bca \frac{(a_{1i})_{x_j}}{a_{1j}}, & 1\leq  i\not= j,\\ 0, & i=j,\eca \quad \d=\diag(\rd x_1, \ldots, \rd x_k),$$
and $\ba_i=(a_{i1}, \ldots, a_{ik})$ for all $1\leq i\leq k$.  Then
\ben
\item fundamental forms of $Y_i$ are 
$$\I_i=\sum_{i=1}^k a_{ij}^2 \rd x_j^2, \quad \II_i= \sum_{m, j=1}^{k, n-k} a_{im} h_{mj} e_{k+j}$$
for some $\calM_{k, n-k}$ matrix $h=(h_{ij})$ (so  $w_{i, k+j}= h_{ij} \rd x_i$ for $1\leq i\leq k$ and $1\leq j\leq n-k$),
\item $w= (w_{ij})_{1\leq i,j\leq n}:=g^{-1}\rd g = \bpm \d F- F^t \d & \d h\\ -h^t \d &0 \epm$ is flat, 
\item  $\rd \ba_i = \ba_i (\d F^t - JF\d J)$ for $1\leq i\leq k$, where $J=\I_{k-\ell, \ell}$, in other words, $\ba_i$ is a parallel field for the $o(k-\ell, \ell)$-valued connection $1$-form $\d F^t- JF\d J$. 
\item $\tau:=\d F^t - JF\d J$ is a flat $o(k-\ell, \ell)$ connection $1$-form, 
\item 
\beq\label{ey}
\o_\l= \bpm w& -\l D^tJ\\ \l D & \tau\epm
\eeq
 is flat for all $\l\in \C$, where $D=(\d, 0)$ and $\d=\diag(\rd x_1, \ldots, \rd x_k)$,
\item let $\bpm g_1 & 0 \\ 0 & g_2\epm$ be a frame of $\o_0=\bpm w & 0\\ 0 & \tau\epm$, then there exists a constant $C\in GL(k)$ such that 
$\rd Y= g_1 Dg_2^{-1} C$. 
\een
\ethm

The equation for $k$-tuples in $\R^n$ of type $\R^{k-\ell, \ell}$ is the equation for $(F, h): \R^k\to gl(k)_\ast\times \calM_{k\times (n-k)}$ such that $w$ and $\tau$ defined in Theorem \ref{eg}  are flat, i.e.,
\beq\label{ex}
\bca \rd w+ w\wedge w=0, & w= \bpm \d F- F^t \d &\d h\\ -h^t\d &0\epm,\\
\rd \tau+ \tau\wedge \tau=0, & \tau=\d F^t - JF\d J, 
\eca
\eeq
where 
$$gl(k)_\ast=\{(y_{ij})\in gl(k)\n y_{ii}=0 \, \, \forall\,\, 1\leq i\leq k\} $$ 
and $J= \I_{k-\ell, \ell}$.  
So $\o_\l$ defined by \eqref{ey} is the Lax pair for the equation \eqref{ex} of $k$-tuples in $\R^n$ of type $\R^{k-\ell, \ell}$.  

\ss\ni {\it Construction of $k$-tuples from parallel frames}

Let $(F,h)$ be a solution of \eqref{ex}, and $E$ is the normalized parallel frame for $\o_\l$ defined by  \eqref{ey}.  Since $\o_0=\bpm w&0\\ 0&\tau\epm$,  $E(x,0)= \bpm g_1 & 0\\ 0 & g_2\epm$ for some $g_1\in O(n)$ and $g_2\in O(k-\ell, \ell)$. Similar argument as for Christoffel pairs of isothermic surfaces gives 
\ben
\item $g_1 D g_2^{-1}$ is closed, so there exists $Y$ such that $\rd Y= g_1Dg_2^{-1}$, 
\item $YC$ is a $k$-tuple in $\R^n$ of type $\R^{k-\ell, \ell}$ for a constant $C\in GL(k)$. 
\item $\frac{\p E}{\p \l} E^{-1}\big|_{\l=0}=\bpm 0& Z\\ -JZ^t &0\epm$ for some $3\times 2$ valued map $Z$ and $Y= Z+c_0$ for some constant $c_0\in \R^n$.
\een

\bs
\section{Soliton hierarchies constructed from symmetric spaces} 

We review the method for constructing soliton hierarchies from a splitting of a Lie algebra (cf. \cite{TU09a}).  

\bdefn
Let $L$ be a formal Lie group, $\calL$ its Lie algebra, and $L_\pm$ subgroups of $L$ with Lie subalgebras $\calL_\pm$.  The pair $(\calL_+, \calL_-)$ is called a {\it splitting\/} of $\calL$ if
$\calL= \calL_+ \oplus \calL_-$ as a direct sum of linear subspaces and $L_+\cap L_-=\{e\}$, where $e$ is the identity in $L$. We call the set  $\co=(L_+L_-)\cap (L_-L_+)$ the {\it big cell\/} of $L$. In other words, $f\in \co$ if and only if $f$ can be factored uniquely as $f_+f_-$ and $g_-g_+$ with $f_\pm, g_\pm\in L_\pm$.  
\edefn

\bthm\label{bv} {\bf (Local Factorization Theorem)} \cite{PreSeg86, TU09a}\par

Suppose $L$ is a closed subgroup of the group of Sobolev $H^1$- loops in a finite dimensional Lie group $G$, and $(\calL_+, \calL_-)$ is a splitting of the Lie algebra $\calL$.  Let $\calV$ be an open subset in $\R^N$, and $g:\calV\to L$ a map such that $(x, \l)\mapsto g(x)(\l)$ is smooth.  If $p_0\in \calV$ and $g(p_0)= k_+k_-= h_-h_+$ with $k_\pm, h_\pm\in L_\pm$, then there exist an open subset $\calV_0\subset \calV$ containing $p_0$ and unique $f_\pm, g_\pm : \calV_0\to L_\pm$ such that $g= g_+g_-= f_-f_+$ on $\calV_0$ and $g_\pm(p_0)= k_\pm$, $f_\pm(p_0)= h_\pm$.
\ethm

\bdefn  A commuting sequence $\calJ=\{J_i\n i\geq 1, \, {\rm integer}\}$ in $\calL_+$ is called a {\it vacuum sequence\/} of the splitting $(\calL_+, \calL_-)$ if
$\calJ$ is linearly independent and each $J_j$ is an analytic function of $J_1$.  
\edefn  

\ss\ni{\bf Construction of soliton hierarchy}

Let $(\calL_+, \calL_-)$ be a splitting of $\calL$, and $\{J_j\n j\geq 1\}$ a vacuum sequence. For  $\xi\in \calL$, let $\xi_\pm$ denote the projection of $\xi$ onto $\calL_\pm$ with respect to $\calL= \calL_++\calL_-$.  Set
 \beq\label{ca}
 \calM= \{ (g^{-1}J_1 g)_+\n g\in L_-\}.
 \eeq   Assume that given smooth $\xi:\R\to \calM$, there is a unique $Q_j(\xi)\in \calL$ such that 
\ben 
\item $[\p_x+\xi, Q_j(\xi)]=0$,
\item $Q_j(\xi)$ is a function of $\xi$ and the derivatives of $\xi$,
\item $Q_j(\xi)$ is conjugate to $J_j$ and $Q_j(J_1)= J_j$.  
\een
Claim that
\beq\label{ax}
\frac{\p \xi}{\p t_j}= [\p_{x} +\xi, (Q_j(\xi))_+]
\eeq
is a PDE system on $\calM$. We only need to show that the right hand side is tangent to $\calM$ at $\xi$: Since $[\p_x+\xi, Q_j]=0$, the right hand side of \eqref{ax} is equal to $-[\p_x + \xi, (Q_j)_-]$.  But it should be in $\calL_+$, so it is equal to $-[\xi, (Q_j)_-]_+$, which is tangent to $\calM$.  Hence this defines a flow on $\calM$.  
 We call \eqref{ax}
the {\it $j$-th flow} and the collection of these flows
the {\it soliton hierarchy constructed from $(\calL_+, \calL_-)$ and $\{J_j\n j\geq 1\}$}.  

\bprop\label{be}
The following statements are equivalent for $\xi: \R^2\to \calM$:
\ben 
\item[(1)] $\xi$ is a solution of the flow \eqref{ax},
\item[(2)] $[\p_x + \xi, \, \p_{t_j} + (Q_j(\xi))_+]=0$,
\item[(3)] $\xi \rd x+ (Q_j(\xi))_+ \rd t_j$ is a flat $\calL_+$-valued connection $1$-form. 
\een
So (3) is a Lax pair of the flow \eqref{ax}.
\eprop

 If $\calL$ is a Lie subalgebra of the Lie algebra of formal power series $A(\l)=\sum_{i\geq n_0} A_i \l^i$ with $A_i\in \calG$ a finite dimensional simple Lie algebra, then equation \eqref{ax} is a PDE with a parameter $\l$.  For examples given in this article, it follows from $[\p_x+ \xi, Q_j(\xi)]=0$ that \eqref{ax} gives a determined PDE system in $\xi$.

\ms\ni {\bf Commuting flows on $L_-$}

Given a splitting $(\calL_+, \calL_-)$ of $\calL$ and a vacuum sequence $\{J_j\n j\geq 1\}$, we consider a hierarchy of flows on the negative group $L_-$:
\beq\label{au}
\frac{\p M}{\p t_j}= -M(M^{-1}J_j M)_-.
\eeq 
A direct computation implies that \eqref{au} are commuting flows on $L_-$, i.e., 
$$-(P_j)_{t_k} + (P_k)_{t_j}+ [P_j, P_k]=0$$ for all $j, k$, where $P_j= -(M^{-1}J_j M)_-$.  
Use $[J_1, J_j]=0$ and a straight forward computation to get the following known results (cf. \cite{TU09a}) :

\bprop \label{cm} If $M(t_1, t_j)$ solves the first and the $j$-th flows \eqref{au} on $L_-$, then 
$(M^{-1} J_1M)_+$ is a solution of the $j$-th flow \eqref{ax}.
\eprop

\bthm
The flows in the soliton hierarchy constructed from a splitting and a vacuum sequence commute.
\ethm

\ss\ni {\bf Formal inverse scattering} \cite{TU00a}

Given an element $f\in L_-$, we use the Local Factorization Theorem to construct a solution of the flow in the soliton hierarchy generated by $J_j$ as follow:  First note that $J_1$ is in the phase space $\calM$ defined by \eqref{ca} and \eqref{ax} is satisfied, i.e., the constant map $J_1$ is the solution of all flows in the hierarchy.  The Lax pair of the flow generated by $J_j$ is $J_1 \rd x+ J_j\rd t_j$.  
Let $E(x,t_j)=\exp(xJ_1+ t_j J_j)$, i.e., $E$ is the normalized parallel frame of the solution on $L_+$ satisfying
$$E^{-1}E_x= J_1, \quad E^{-1} E_{t_j} = J_j, \quad E(0, \l)= \I.$$
By Theorem \ref{bv}, given $f\in L_-$, we can factor
$$f^{-1}E(x, t_j)= \ti E(x, t_j) \ti f(x,t_j)^{-1}$$
with $\ti E(x, t_j)\in L_+$ and $\ti f(x, t_j)\in L_-$ for $(x,t_j)$ in some open subset of the origin.  We claim that $\ti f$ is a solution of \eqref{au} for the first and the $j$-th flow.  To see this, note that $\ti E= f^{-1} E \ti f$ and
$$\ti E^{-1}\ti E_x= \ti f^{-1} J_1\ti f + \ti f^{-1} \ti f_x, \quad \ti E^{-1}\ti E_{t_j}= \ti f^{-1}J_j \ti f + \ti f^{-1} \ti f_{t_j}.$$
Since the left hand sides are in $\calL_+$ and $\ti f^{-1}\ti f_x, \ti f^{-1}\ti f_{t_j}$ are in $\calL_-$, the above equation implies that 
$$\bca \ti f^{-1}\ti f_x= -(\ti f^{-1}J_1 \ti f)_-, &\\ \ti f^{-1}\ti f_{t_j}= -(\ti f^{-1} J_j \ti f)_-.\eca$$
Hence $\ti f(x, t_j)$ is a solution of \eqref{au} and this proves the claim.  By Proposition \ref{cm}, $\xi = (\ti f^{-1}J_1 \ti f)_+$ is a solution of the flow generated by $J_j$.

\beg \label{bf} {\bf The $G$-hierarchy\/} \cite{AKNS74, Sa84, Wil91,TU00a}\hfil\par

Let $G$ be a complex simple Lie group, and $L(G)$ the group of  smooth loops $f:S^1\to G$, $L_+(G)$ the subgroup of $f\in L(G)$ that can be extended holomorphically to $|\l|<1$, and $L_-(G)$ the subgroup of $f\in L(G)$ that can be extended holomorphically to $\infty=|\l|>1$ and $f(\infty)=\I$. 
   The corresponding Lie algebras are 
\begin{align*}
\calL(\calG) &= \{A(\l)=\sum_i A_i \l^i\n A_i\in \calG\},\\
\calL_+(\calG)&=\{A\in \calL(\calG)\n A(\l)=\sum_{j\geq 0} A_j \l^j\},\\
\calL_-(\calG)&=\{A\in \calL(\calG)\n A(\l)=\sum_{j<0}A_j \l^j\},
\end{align*}
where $\calG$ is the Lie algebra of $G$.  
Then  $(\calL_+(\calG), \calL_-(\calG))$ is a splitting of $\calL(\calG)$

Let $\calA$ be a maximal abelian subalgebra of $\calG$, and  $\calA^\perp$ the orthogonal complement of $\calA$ with respect to the Killing form $(\, , )$ of $\calG$.
The dimension of $\calA$ is the {\it rank\/} of $G$. An element $\xi\in \calG$ is {\it regular\/} if $\ad(\xi)$ is semi-simple and the centralizer $\calG_\xi$ is a maximal abelian subalgebra. If $\xi$ is regular, then $\ad(\xi)$ is a linear isomorphism of $\calG_\xi^\perp$. Let  $\{a_1, \ldots, a_{r}\}$ be a basis of $\calA$ such that  $a_1$ is regular.  Then
$$\calJ=\{J_{i,j}= a_i\l^j\n 1\leq i\leq r, j\geq 1\}$$
is a vacuum sequence with $J_1= J_{1,1}= a_1\l$. A direct computation shows that $\calM$ defined by \eqref{ca} is 
$$\calM=J_1+([a_1\l, \calL_-])_+=J_1+ \calA^\perp.$$

To write down the flow generated by $J_{i,j}$, we construct 
$$Q_i(u)= a_i \l +\sum_{k\leq 0}  Q_{i,k}(u)\l^k$$
satisfying
\beq\label{bh}
\bca [\p_x + a_1\l + u, Q_{i}(u)]=0,\\ f_j(Q_{i}(u))=f_j(a_i\l), & 1\leq j\leq r,\eca
\eeq
where $f_1, \ldots, f_r$ are free generators of the ring of invariant polynomials on $\calG$ (for example, if $\calG= sl(n)$, then $r=n-1$ and $f_j(A)$ can be chosen to be $\tr (A^j)$ for $2\leq j\leq n$). 
Equate the coefficient of $\l^k$ in the first equation of \eqref{bh} to get the recursive formula
\beq\label{bk}
(Q_{i,k})_x+ [u, Q_{i, k}] + [a_1, Q_{i, k-1}]=0.
\eeq
We use \eqref{bk} and the second equation of \eqref{bh} to prove that $Q_{i, k}$ is a polynomial differential operator of $u$.  Since $M^{-1}J_{i,j} M= \l^{j-1}M^{-1}J_{i,1} M$, the flow generated by $J_{i,j}$ is \eqref{ax}, i.e., 
$$\frac{\p (a_1\l + u)}{\p t_{i,j}}= [\p_x+a_1\l + u, \, a_i\l^j + Q_{i, 0}\l^{j-1} + \cdots + Q_{i, 1-j}].$$
Although the right hand side is a degree $j+1$ polynomial  in $\l$, it follows from the recursive formula \eqref{bk} that all coefficients of $\l^k$ of the right hand side are zero except the constant term. So the flow equation generated by $J_{i,j}$ is the following PDE for $u$:
\beq\label{bz}
u_{t_{i,j}}= [\p_x + u, Q_{i, 1-j}]= [Q_{i, -j}, a_1].
\eeq
By Proposition \ref{be},  equation \eqref{bz} has a Lax pair
$$\o_\l = (a_1\l + u) \rd x + (a_i\l^j + Q_{i, 0}\l^{j-1} + \cdots + Q_{i, 1-j})\rd t_{i,j}.$$  
We call this hierarchy of flows the {\it $\calG$-hierarchy}. For example, for general $\calG$, the flow generated by $J_{1,j}$ in the $\calG$-hierarchy is the PDE for $u:\R^2\to \calA^\perp$: 
$$u_{t_{1,j}}= \ad(a_j)\ad(a_1)^{-1}(u_x) + [u, \ad(a_j)\ad(a_1)^{-1}(u)],$$ 
 and its Lax pair is 
 $$\o_\l = (a_1\l + u) \rd x + (a_j\l + \ad(a_j)\ad(a_1)^{-1}(u))\rd t_{1,j}.$$
\eeg

\beg {\bf The $U$-hierarchy} \cite{TU00a}\par

Let $\tau$ be a Lie group involution of $G$ such that $\rd \tau_e:\calG\to \calG$  
(still denoted by $\tau$) is conjugate linear.  Let $U$ denote the fixed point set of $\tau$, and $\calU$ the Lie algebra of $U$, i.e., $\calU$ is a {\it real form\/} of $\calG$.  Let $L^\tau(G)$ denote the subgroup of all $f\in L(G)$ satisfying the {\it $U$-reality condition\/}
 \beq\label{an}
\tau(f(\bar\l))=f(\l),
\eeq
and $L^\tau_\pm(G)= L^\tau(G)\cap L_\pm(G)$. Let $\calL^\tau(\calG)$ and $\calL_\pm^\tau(\calG)$ denote the corresponding Lie algebras.  Let $\{a_1, \ldots, a_n\}$ be a basis of a maximal abelian subalgebra of $\calU$ such that $a_1$ is regular, and $\calJ=\{J_{ij}=a_i \l^j\n 1\leq i\leq n, j\geq 1\}$. Then $(\calL_+^\tau(\calG), \calL_-^\tau(\calG))$ is a splitting and $\calJ$ is a vacuum sequence.  The flows generated by $J_{ij}$'s form the {\it $U$-hierarchy} and flows in the $U$-hierarchy are evolution equations on $C^\infty(\R, \calA^\perp\cap\calU)$.  
For example, for  $\tau(g) = (\bar g^t)^{-1}$ on $sl(2,\C)$. Then  $\calU= su(2)$. Let $a= \diag(\im, -\im)$. The flows are evolution PDE on  $C^\infty(\R, Y)$, where $Y=\left\{\bpm 0 & q\\ -\bar q & 0\epm\n q\in \C\right\}$ and  the flow generated by $J_{1,2}= a\l^2$ in
the $su(2)$-hierarchy is the NLS 
\eeg

\beg \label{bi} {\bf The $\frac{U}{K}$-hierarchy} \cite{TU00a}\par

Let $\tau, \sigma$ be commuting involutions of $G$ such that the induced involutions $\tau$ and $\sigma$ on $\calG$ are conjugate and complex linear  respectively, and $U$ the fixed point set of $\tau$ on $G$ and $K$ the fixed point set of $\sigma$ on $U$ (so $\frac{U}{K}$ is a symmetric space).  Let $\calP$ denote the $-1$ eigenspace of $\sigma$ in $\calU$. Then we have $\calU=\calK+\calP$ and 
$$[\calK, \calK]\subset \calK, \quad [\calK, \calP]\subset \calP, \quad [\calP, \calP]\subset \calK.$$
This is the {\it Cartan decomposition\/} for $\frac{U}{K}$. Note that $K$ acts on $\calP$ by conjugation.  
An element $b\in \calP$ is {\it regular\/} if the $K$-orbit of $b$ in $\calP$ is maximal.  If $b$ is regular, then $\{\xi\in \calP\n [b,\xi]=0\}$ is a maximal abelian subalgebra and is the kernel of $\ad(b):\calP \to \calK$.
 
Let $\calA$ be a maximal abelian subalgebra in $\calP$, and $\{ a_1, \ldots, a_n\}$ a basis of $\calA$ such that $a_1$ is regular (i.e., $\ad(a_1)$ is a linear isomorphism from $\calK\cap \calK_{a_1}^\perp$ onto $\calP\cap \calA^\perp$, where $\calK_{a_1}= \{k\in K\n [a_1, k]=0\}$.  The dimension of $\calA$ is the {\it rank\/} of the symmetric space.  

Let $\calL^{\tau, \sigma}(\calG)$ be the subalgebra of $\xi\in L(\calG)$ satisfying the {\it $\frac{U}{K}$-reality condition\/}
\beq\label{ao}
\tau(\xi(\bar\l))= \xi(\l), \qquad \sigma(\xi(-\l))= \xi(\l),
\eeq
and 
$$\calL_\pm^{\tau,\sigma}(\calG)= \calL^{\tau,\sigma}(\calG)\cap \calL_\pm(\calG)).$$
Then $(\calL^{\tau,\sigma}_+(\calG), \calL_-^{\tau,\sigma}(\calG))$ is a splitting and 
$$\calJ=\{ J_{ij}= a_i \l^j\n 1\leq i\leq n, j\geq 1 \, \, {\rm odd \,\, integer\/} \}$$
is a vacuum sequence.  The hierarchy constructed from these are called the $\frac{U}{K}$-hierarchy and the flows in this hierarchy are evolution equations on $C^\infty(\R, \calK_{a_1}^\perp)$, where $\calK_{a_1}^\perp= \{ y\in \calK\n (y, k)=0 \, \forall\, \, k\in \calK_{a_1}\}$.  
For example, the symmetric space given by $\tau(g)= (\bar g^{-1})^t$ and $\sigma(g)= (g^t)^{-1}$ on $G=SL(2,\C)$ is $\frac{SU(2)}{SO(2)}=S^2$. Let $a=\diag(\im, -\im)$. The flows in the $\frac{SU(2)}{SO(2)}$-hierarchy are for  $u=\bpm 0& q\\ -q& 0\epm$ and the flow generated by $J_{1,3}= a\l^3$ is the mKdV.
\eeg

\brem
If $\frac{U}{K}$ has {\it maximal rank\/}, i.e., the rank of $\frac{U}{K}$ is equal to the rank of $U$, then:
\ben
\item[\bu] A maximal abelian subalgebra $\calA$ in $\calP$ is also a maximal abelian subalgebra of $\calU$ over $\R$ and is a maximal abelian subalgebra of $\calG$ over $\C$.
\item[\bu] Fix a basis $\{a_1, \ldots, a_n\}$ of $\calA$ over $\R$. The phase space for flows in the $G$-hierarchy is $C^\infty(\R, \calA^\perp)$.
\item[\bu] The  flow generated by $J_{i,j}$ in the $G$-hierarchy leaves  $C^\infty(\R, \calA^\perp\cap \calU)$ invariant and the restricted flows form the $U$-hierarchy.
\item[\bu] The flow generated by $J_{i,2j+1}$ of the $U$-hierarchy leaves the subspace $\calK_{a_1}^\perp$ invariant and the restricted flows form the $\frac{U}{K}$-hierarchy. 
\een 
\erem

\ss\ni {\bf The matrix NLS hierarchy} \cite{ForKul83, TU00a}

Let $\tau(g)= (\bar g^t)^{-1}$ be the involution of $G=GL(n,\C)$ that defines the real form $U=U(n)$, and $(\calL_-^\tau(\calG), \calL_-(\calG))$ the splitting that gives the $U$-hierarchy.  Let $a=\im \I_{k, n-k}$. Then 
$$\calJ=\{a\l^j\n j\geq 1\}$$
is a vacuum sequence. The flows constructed by this splitting and hierarchy are equations for $u:\R^2\to \calM_{k\times (n-k)}$, and the flow generated by $a\l^2$ is the {\it matrix NLS\/}, $q_t= \frac{\im }{2}(q_{xx} + 2 q\bar q^t q)$. 

\ss
\ni {\bf The $-1$ flow associated to $\frac{U}{K}$} \cite{Ter97}

We use the same notation as for the $\frac{U}{K}$-hierarchy.  Given $b\in \calA$, the $-1$ flow associated to $\frac{U}{K}$ is the equation for $g:\R^2\to K$:
\beq\label{cq}
 (g^{-1}g_x)_t = [a, g^{-1} bg].
\eeq
It is easy to check that $g$ is a solution of \eqref{cq} if and only if $\o_\l$ is flat for all $\l\not=0$, where
\beq \label{cs}
\o_\l= (a_1\l + g^{-1}g_x) \rd x + \l^{-1} g^{-1}bg\rd t.
\eeq
For example, the $-1$-flow associated to $\frac{SU(2)}{SO(2)}$ defined by $a= \diag(\im, -\im)$ and $b= -\frac{a}{4}$ is the equation for $g=\bpm \cos q & -\sin q\\ \sin q & \cos q\epm$, \eqref{cs} is \eqref{bb}, and \eqref{cq} gives the SGE.  

\beg \label{ct} {\bf Twisted $\frac{U}{K_1}$-hierarchy} \cite{Ter09}\par

Let $\tau$ be the conjugate involution of the complex simple Lie group $G$ that gives the real form $U$, $\sigma_1$ and $\sigma_2$ involutions of $\calG$ such that $\sigma_1, \sigma_2$ and $\tau$ commute, and 
$$\calU= \calK_1+ \calP_1, \quad \calU=\calK_2+\calP_2$$
Cartan decompositions for $\sigma_1$ and $\sigma_2$ respectively.  Let $\calA$ be a maximal abelian subalgebra in $\calP_1$.  Assume that 
\ben
\item[1.] $\sigma_2(\calA)\subset \calA$,
\item[2.] $K_1\cap K_2= S_1\times S_2$, $K_1= S_1\times K_1'$, $K_2= K_2'\times S_2$ as direct product of subgroups. 
\een
Let $L=L^{\tau,\sigma_1}$ denote the group of holomorphic maps $f$ from $\e<|\l | <\e^{-1}$ to $G$ satisfying the $U/K_1$-reality condition: 
$$\tau(f(\bar\l))= f(\l), \quad \sigma_1(f(-\l))= f(\l).$$
Let $L_+$ denote the subgroup of $f\in L$ such that $\sigma_2(f(\l^{-1}))= f(\l)$ and $f(1)\in K_2'$, and $L_-$ the subgroup of $f\in L$ that can be extended holomorphically to $\infty\geq |\l | > \e$ and $f(\infty)\in K_1'$.  Then $L_+\cap L_-=\{e\}$ and the Lie algebras are:
\begin{align*}
&\calL=\{ \xi(\l)=\sum_j \xi_j \l^j\,\big|\, \xi_j\in \calK_1 \, {\rm if\, } k {\rm \, is \, even,}\, \in \calP_1, \, {\rm if\, } k \, {\rm is\, odd. }\},\\
&\calL_+=\{ \xi\in \calL\n \xi_{-j}= \sigma_2(\xi_j), \xi(1)\in \calK_2'\},\\
&\calL_-=\{\xi(\l)=\sum_{j\leq 0} \xi_j\l^j\in \calL\n \xi_0\in \calK_1'\}.
\end{align*}
Let $\{a_1, \ldots, a_n\}$ be a basis of $\calA$ such that $a_1$ is regular with respect to the Ad$(K_1)$ action on $\calP_1$, and $\calJ=\{J_{i,j}\n 1\leq i\leq n, j\geq 1 \, {\rm odd\/}\},$ where
$$J_{i, j} = a_i \l^j + \sigma_2(a_i)\l^{-j}.$$
Then $(\calL_+, \calL_-)$ is a splitting of $\calL$ and $\calJ$ is a vacuum sequence.  We call the hierarchy constructed from this splitting and vacuum sequence {\it the $\frac{U}{K_1}$-hierarchy twisted by $\sigma_2$\/}. 
The phase space of this hierarchy is $C^\infty(\R, \calM)$, where 
$$\calM= \{g^{-1}a_1g\l + v+\sigma_2(g^{-1}a_1g)\l^{-1}\n g\in K_1', v\in \calS_1\}.$$
\eeg

\beg\label{br}
{\bf A twisted $\frac{O(n,n)}{O(n)\times O(n)}$-hierarchy} \cite{Ter09}\par

 Let $\calG= o(n, n, \C)$, $\tau(x)= \bar x$, and
 $$\sigma_1(x)= \I_{n,n} x \I_{n,n}^{-1}, \quad \sigma_2(x)= \I_{n+1, n-1} x \I_{n+1, n-1}^{-1}.$$ Then 
\begin{align*}
&\calU= o(n,n), \quad \calK_1= o(n)\times o(n), \quad \calK_2= o(n,1)\times o(n-1),\\
& \calK_1\cap \calK_2= \calS_1+\calS_2, \, {\rm where }\quad \calS_1=o(n)\times 0, \, \calS_2= 0\times o(n-1),\\
&\calK_2'= \calS_1 +(\calP_1\cap \calK_2)= o(n,1),\quad \calK_1'= 0\times o(n).
\end{align*}
The space 
$$\calA=\left\{\bpm 0& D\\ D&0\epm\, \bigg|\, D\in gl(n,\R)\,\, {\rm is\,\, diagonal\/}\right\}$$
is a maximal abelian subalgebra in $\calP_1$ and $\sigma_2(\calA)\subset \calA$.  Choose a basis $\{a_1, \ldots, a_n\}$ of $\calA$ such that $a_1$ is regular. Then $\tau, \sigma_1, \sigma_2$ satisfy all the conditions given above and we obtain the $\frac{O(n,n)}{O(n)\times O(n)}$-hierarchy twisted by $\sigma_2$. 
\eeg

\ms
Next we give a brief discussion of bi-Hamiltonian structure, conservation laws, and formal inverse scattering for the $U$-hierarchy. 

\ss\ni {\bf Bi-Hamiltonian structure for the $U$-hierarchy} (cf. \cite{DS84, Ter97})

Let $(\, , )$ denote a bi-invariant non-degenerate bilinear form on $\calU$, and 
$$\li u, v\ri = \int_{-\infty}^\infty (u(x)v(x))\rd x$$
the induced bi-linear form on $V=\calS^\infty(\R, \calA^\perp)$ the space of Schwartz maps from $\R$ to $\calA^\perp$. 
Given a functional $F$ on $V$,  the gradient of $F$ is defined by 
$$dF_u(v)= \li \K F(u), v\ri,$$
(i.e., $\K F(u)=0$ is the Euler-Lagrangian equation for $F$).  A Poisson structure on $V$ is an operator $J: V\to L(V,V)$, $u\mapsto J_u$ such that 
$$\{F_1, F_2\}(u) = (J_u(\K F_1(u)), \K F_2(u))$$
defines a Lie bracket on $V$ and $\{\, \}$ satisfies the product rule.  The Hamiltonian equation for $F:V\to \R$ is 
$$\frac{\rd u}{\rd t}= J_u(\K F(u)).$$

Two Poisson structures $\{\, , \, \}_1, \{\, , \, \}_2$ on $V$ are {\it compatible\/} if 
$$c_1 \{\, , \, \}_1 + c_2 \{\, , \, \}_2$$
is again a Poisson structure for any constant $c_1, c_2$.  
 
 Given a smooth map $u:\R\to \calA^\perp$, let $P_u$ be the operator on $C(\R, \calA^\perp)$ defined by 
$$P_u(v)= (\ti v)_x + [u,\ti v], \qquad \ti v= v+ T, \quad T(x)= -\int_{\infty}^x [u, v]_0,$$
where $\xi_0$ and $\xi^\perp$ denote the projection of $\xi\in \calU$ to $\calA$ and $\calA^\perp$ respectively.  By definition, $P_u(v)\in \calA^\perp$.  
Let $J_0$ and $J_1$ be the operator from $V$ to $L(V,V)$ defined by  
$$(J_0)_u= -\ad(a_1), \quad (J_1)_u = P_u.$$
 Define
$$\{F, G\}_0(u)= \li [\K F(u), a_1], \K G(u)\ri,\quad
\{F, G\}_1(u)= \li [P_u(\K F(u)), \K G(u)\ri.$$
The following are known (cf. \cite{DS84}, \cite{Ter97}):
\ben
\item $\{\, , \, \}_0$ and $\{\, , \, \}_1$ are compatible Poisson structures on $C(\R, \calA^\perp)$. 
\item Set
\beq\label{bha} F_{i,j}(u)=-\frac{1}{j} \int_{-\infty}^\infty (Q_{i,-j}(u), a_1)\rd x.
\eeq
Then $\K F_{i,j}(u)= Q_{i, -j+1}(u)^\perp$ and the flow generated by $J_{i, j}$ is 
$$u_t= J_0(\K F_{i, j+1}) = J_1(\K F_{i, j}).$$
\item Both Poisson structures can be constructed from the natural Poisson structures of co-adjoint orbits of $L_-^\tau(G)$. 
\een

\bs
\section{The $\frac{U}{K}$-system and the Gauss-Codazzi equations}\par

We review the definition of the $\frac{U}{K}$-system, the twisted $\frac{U}{K}$-system, and the $-1$ flow on the $\frac{U}{K}$-system and see that SGE, GSGE, equations for isothermic surfaces, for $k$-tuples in $\R^n$ of type $\R^{k-\ell, \ell}$, and for flat Lagrangian submanifolds in $\C^n$ are $\frac{U}{K}$-systems. 

\ss\ni {\bf The $\frac{U}{K}$-system} \cite{Ter97}\par\hfil

Let $\frac{U}{K}$ be a rank $n$ symmetric space, $\calU= \calK+ \calP$ a Cartan decomposition, $\calA$ a maximal abelian subspace in $\calP$, and $\{a_1, \ldots, a_n\}$ a basis of $\calA$.  The {\it $\frac{U}{K}$-system} is the following over-determined first order non-linear PDE system for $v:\R^n\to \calA^\perp\cap \calP$:
\beq \label{ap}
[a_i, v_{t_{1,j}}] -[a_j , v_{t_{1,i}}]=[[a_i, v], [a_j, v]], \qquad 1\leq i\not= j\leq n,
\eeq

It follows from the definition that the following statements are equivalent for $v:\R^n\to \calA^\perp\cap \calP$:
\ben 
\item $v$ is a solution of \eqref{ap},
\item the following connection $1$-form on $\R^n$ is flat for all parameters $\l\in \C$:
\beq \label{aq}
\o_\l = \sum_{i=1}^n (a_i \l + [a_i, v]) \rd x_i
\eeq
($\o_\l$ is a Lax pair of the $\frac{U}{K}$-system),
\item $\o_s$ is flat for some $s\in \R\cup i\R$,
\item if $a_1$ is regular, then $u=[a_1,v]$ is a solution of the flow generated by $J_{i, 1}= a_i \l$ in the $\frac{U}{K}$-hierarchy. 
\een

\brem If we use a different basis of $\calA$, the $\frac{U}{K}$-systems differ by a linear coordinate change.  If two maximal abelian subalgebras $\calA$ and $\ti\calA$ are conjugated by an element in $K$, then the corresponding $\frac{U}{K}$-systems are equivalent.  If $\frac{U}{K}$ is a Riemannian symmetric space, then any two maximal abelian subalgebras in $\calP$ are conjugate by an element of $K$, so there is a unique $U/K$-system. But when $\frac{U}{K}$ is a pseudo-Riemannian symmetric space, there may be more than one maximal abelian subalgebras in $\calP$ modulo the conjugation action of $K$ on $\calP$. Hence there may be more than one non-equivalent $\frac{U}{K}$-system associated to $\frac{U}{K}$.  
\erem 

Statement (4) given above means that the $\frac{U}{K}$-system combines the commuting flows in the $\frac{U}{K}$-hierarchy generated by $J_{1, 1}=a_1\l, \ldots, J_{n,1}=a_n\l$ together.  

\ss\ni{\bf Curved flats in symmetric spaces} \hfil \par

Recall that a {\it flat\/} of a symmetric space $\frac{U}{K}$ is a totally geodesic flat submanifold of $\frac{U}{K}$.  If $\calA$ is a maximal abelian subalgebra in $\calP$, then $A=\exp(\calA)K$ is a flat through $eK$ and $gA$ is a flat through $gK$. Moreover, all flats are obtained this way. 

\bdefn\cite{FP96a}
A {\it curved flat\/} in $\frac{U}{K}$ is an immersed flat submanifold of $\frac{U}{K}$ that is tangent to a flat of $\frac{U}{K}$ at every point.  
\edefn

\bdefn\cite{Ter08}
Let $\frac{U}{K}$ be a symmetric space, and $\calU=\calK+\calP$ a Cartan decomposition.  A flat submanifold $M$ of $\calP$ is called an {\it abelian flat submanifold\/} if $TM_x$ is a maximal abelian subalgebra of $\calP$ for all $x\in M$. Here the metric on $\calP$ is the restriction of the Killing form $(\, ,)$ of $\calU$ to $\calP$.
\edefn

If we identify the tangent space of $\frac{U}{K}$ at $eK$ to be $\calP$, then a flat submanifold $\Sigma$ in $\frac{U}{K}$ is a curved flat if and only if $g^{-1}T\Sigma_{gK}$ is a maximal abelian subalgebra of $\calP$ for all $gK\in \Sigma$.  A curved flat $\Sigma$ is semi-simple if $g^{-1}T\Sigma_{gK}$ is a semi-simple maximal abelian subalgebra of $\calP$ for all $gK\in \Sigma$. 

Let $\frac{U}{K}$ be the symmetric space defined by $\tau, \sigma$. Then the map $\frac{U}{K}\to U$ defined by $gK\mapsto g\sigma(g)^{-1}$ is well-defined and gives an isometric embedding of the symmetric space $\frac{U}{K}$ into $U$ as a totally geodesic submanifold.  This is called the {\it Cartan embedding\/} of $\frac{U}{K}$ in $U$.

The following is known (\cite{FP96a, Ter08}): 

\bthm\label{fh}
 Suppose $v$ is a solution of the $\frac{U}{K}$-system and $E$ is its a parallel frame. Then:
 \ben
 \item 
 $Y= E(x,\l)\sigma(E(x, \l))^{-1}\, \bigg|_{\l=1} = E(x, 1)E(x, -1)^{-1}$
is a curved flat. Conversely, all local semi-simple curved flats can be constructed this way.  In other words, {\it the $\frac{U}{K}$-system can be viewed as the equation for curved flats in $\frac{U}{K}$ with a ``good coordinate system''}.   
\item $Z= \frac{\p E}{ \p \l}E^{-1} \big|\ _{\l=0}$  is an abelian flat in $\calP$.  Conversely, locally all abelian flats in $\calP_0$ can be constructed
this way, where $\calP_0$ is the subset of regular points in $\calP_0$.  
\een
\ethm

\beg{\bf The $\frac{U(n)\sdp \C^n}{O(n)\times \R^n}$-system} \cite{TW08}

Let $U(n)\sdp \C^n$ denote the group of unitary rigid motions of $\C^n=\R^{2n}$, and $\calG$ the complexified $u(n)\sdp \C^n$, i.e., 
$$\calG=\left\{\bpm b&c&x\\-c & b& y\\ 0&0&0\epm\, \bigg|\, b^t=-b, c^t=-c, b, c\in gl(n,\C), x, y\in \C^n\right\}.$$
Let $\tau, \sigma:G\to G$ be the involutions defined by
$$\tau(g)=\bar g, \quad \sigma(g)= TgT^{-1}, \quad {\rm where\, } T=\bpm \I_n &0&0\\ 0& -\I_n &0\\ 0&0&1\epm.$$
The fixed point set of $\tau$ is $U(n)\sdp \C^n$, $\sigma$ and $\tau$ commute, and the corresponding symmetric space is $\frac{U(n)\sdp \C^n}{O(n)\sdp \R^n}$.  The Cartan decomposition is $u(n)\sdp \C^n = \calK +\calP$, where
\begin{align*}\calK&=\left\{\bpm b&0&x\\0 & b& 0\\ 0&0&0\epm\,\bigg|\, b\in o(n),  x\in \R^n\right\}, \\
\calP&=\left\{\bpm 0&-c&0\\ c & 0& y\\ 0&0&0\epm\,\bigg|\, c=c^t, \bar c=c, y\in \R^n  \right\}.
\end{align*}
Then $\{a_i=e_{n+i, i} -e_{i, n+i}\n 1\leq i\leq n\}$ form a basis of a maximal abelian algebra $\calA$ in $\calP$.  The $\frac{U(n)\sdp \C^n}{O(n)\sdp \R^n}$-system is the system for 
$q=\bpm 0& \b &0\\ -\b & 0 & -h\\ 0&0&0\epm$ 
 given by the condition that 
\beq\label{du}
\o_\l= \sum_{i=1}^n(a_i \l + [a_i, q])\rd x_i=\bpm [\d, \b]& \l \d & \d h\\ -\l \d & [\d,\b] &0\\ 0&0&0\epm
\eeq
is flat for all $\l\in \C$.  Note that this is the Lax pair \eqref{dm} for flat Lagrangian submanifolds in $\C^n$.
\eeg

\beg {\bf The $\frac{O(4,1)}{O(3)\times O(1,1)}$-system}\par

The involutions that gives $\frac{O(4,1)}{O(3)\times O(1,1)}$ is $\tau(g)= \bar g$ and $\sigma(g)= \I_{3,2}g\I_{3,2}^{-1}$, and the Cartan decomposition is $o(4,1)= \calK+\calP$ with $\calK= o(3)\times o(1,1)$ and 
$$\calP=\left\{\bpm 0& \xi\\ -J\xi^t &0\epm\,\big|\, \xi \,\, {\rm is\, a\, real \, } 3\times 2 \, {\rm matrix\/}\right\}, \quad J= \diag(1,-1).$$  
Note that 
$$\calA=\left\{\bpm 0 &\xi \\ -J\xi^t & 0\epm\,\bigg|\, \xi =\bpm c_1&0\\ 0& c_2 \\ 0&0\epm\right\}$$ 
is a maximal abelian subalgebra in $\calP$.  Let
 $\{a_1, a_2\}$ be a basis of $\calA$ defined by
$$a_i= \bpm 0& D_i\\ -JD_i^t & 0\epm, \quad D_1=\bpm 1&0\\ 0&0\\0&0\epm, \quad D_2= \bpm 0&0\\0&1\\0&0\epm.$$
The $\frac{O(4,1)}{O(3)\times O(1,1)}$-system  \eqref{ap} is for $v=\bpm 0 & \xi \\ -J\xi^t &0\epm$ with $\xi=\bpm 0 & f_1\\ f_2 &0\\ -r_1& r_2\epm$. Write down this system in terms of $f_1, f_2, r_1, r_2$ we get
\beq\label{cb}\bca (f_1)_{x_2}= -(f_2)_{x_1}, \\ (f_2)_{x_2}-(f_1)_{x_1} - r_1 r_2=0, \\ (r_1)_{x_2}=- f_2 r_2,\\ (r_2)_{x_1}= f_1 r_1.\eca
\eeq
Its Lax pair is 
\beq\label{cc}\o_\l = \bpm -DJ\xi^t +\xi J D^t & D \l\\  -JD^t \l & -JD^t \xi+J\xi^t D\epm, \quad {\rm where \, } D= \bpm \rd x_1& 0\\ 0 & \rd x_2 \\ 0&0\epm. 
\eeq
The first equation of \eqref{cb} implies that there exists $q$ such that $f_1 = q_{x_1}$ and $f_2= -q_{x_2}$. Write \eqref{cb} in terms of $q, r_1, r_2$ we get the Gauss-Codazzi equation \eqref{am} for isothermic surfaces.  Moreover, the Lax pair 
 \eqref{cc} is the Lax pair \eqref{ck} for isothermic surfaces in $\R^3$.  
\eeg

\beg {\bf The $\frac{O(n+k-\ell, \ell)}{O(n)\times O(k-\ell, \ell)}$-system}

We choose 
$$a_i= \bpm 0 & -D_i J\\ D^t &0\epm,  \quad 1\leq i\leq k, \quad v=\bpm 0 & -\xi^tJ\\ \xi & 0\epm,$$
where $D_i^t= (e_{ii}, 0) \in \calM_{k\times n}$, and $e_{ii}$ is the diagonal $k\times k$ matrix with all entries zero except the $ii$-th entry is $1$.  
The $\frac{O(n+k-\ell, \ell)}{O(n)\times O(k-\ell, \ell)}$-system is the PDE for  $\xi=(y, \g):\R^k\to  gl_\ast(k)\times \calM_{k, n-k}$ with Lax pair  $\o_\l=\sum_{i=1}^k (a_i\l + [a_i,v])\rd x_i$. We write $\o_\l$ in terms of $y, \g$ to get
\beq\label{ef}
\o_\l =\bpm w & -\l \eta J \\ \eta^t\l & \tau\epm,  \quad {\rm where}
\eeq 
$$w= \bpm -\d J y+ y^t J \d & \d J\g\\ -\g^t J\d &0\epm, \quad \tau= -\d y^t J+ y\d J, \quad \eta^t= (\d, 0),$$
and $\d= \diag(\rd x_1, \ldots, \rd x_k)$.

Set $F= -Jy$ and $h= J\g$, then \eqref{ef} is the same Lax pair \eqref{ey} for $k$-tuples in $\R^n$ of type $\R^{k-\ell, \ell}$ given in Theorem \ref{eg}.  So the $\frac{O(n+k-\ell, \ell)}{O(n)\times O(k-\ell, \ell)}$-system is the equation for $k$-tuples in $\R^n$ of type $\R^{k-\ell, \ell}$.
\eeg 

\beg\cite{BDPT02} The $\frac{O(5)}{O(3)\times O(2)}$-system is the equation for 

\ben
\item $2$-tuples in $\R^3$ of type $O(2)$,
\item flat surfaces in $S^4$ with flat and non-degenerate normal bundle,
\item surfaces in $S^4$ with constant sectional curvature $1$ and flat and non-degenerate normal bundle. 
\een
Moreover, if $v$ is a solution of the $\frac{O(5)}{O(3)\times O(2)}$-system, and $E$ a parallel frame of the Lax pair of $v$.  Write $E(x,0)= \bpm g_1 & 0\\ 0 & g_2\epm$, and $D=\bpm \d \\ 0\epm$ is $\calM_{3\times 2}$ valued, where $\d= \diag(\rd x_1, \rd x_2)$.  Then: 
\ben
\item $g_1D g_2^{-1}$ is closed, so there exists  $Y=(Y_1, Y_2)\in \calM_{3\times 2}$ such that $\rd Y= g_1 D g_2^{-1}$, and $Y$ is a $2$-tuple of surfaces in $\R^3$ of type $O(2)$.
\item The first column of $E(x,r)\bpm g_1^{-1} &0 \\  0 & \I_2\epm$ is a flat surface in $S^4$ with flat and non-degenerate normal bundle.
\item The third column of $E(x,r)\bpm \I_3 & 0\\ 0 & g_2^{-1}\epm$ is a surface in $S^4$ with constant curvature $1$ and flat, non-degenerate normal bundle.
\een
\eeg

Analogous results hold for $\frac{U}{K}$-system when $\frac{U}{K}$ is a real Grassmannian. 

\beg \cite{TW08}:  The $\frac{U(n)}{O(n)}$-system is the equation for
\ben
\item Egoroff orthogonal coordinate systems of $\R^n$,
\item flat Lagrangian submanifolds of $\C^n$ that lie in $S^{2n-1}$, 
\item flat Lagrangian submanifolds of $\C P^{n-1}$.
\een
\eeg

\ss\ni {\bf Twisted $\frac{U}{K_1}$-system} \cite{FP96b, Bra07b, Ter09}\par

We use the same notation as for twisted $\frac{U}{K_1}$-hierarchy.  The $\frac{U}{K_1}$-system twisted by $\sigma_2$ is the PDE for maps $g:\R^n\to K_1'$ and $v_i:\R^n\to \calS_1$ such that the connection $1$-form 
\beq\label{dn}
\o_\l= \sum_{i=1}^n ((ga_i g^{-1})\l + v_i + \sigma_2(ga_i g^{-1})\l^{-1})\rd x_i
\eeq
 is flat for all non-zero parameters $\l\in \C$.  So the $\frac{U}{K_1}$-system twisted by $\sigma_2$ is given by the collection of flows in the $\frac{U}{K_1}$-hierarchy twisted by $\sigma_2$ generated by $a_i\l + \sigma_2(a_i)\l^{-1}$ for $1\leq i\leq n$.

\ms
\beg\label{bq} {\bf A twisted $\frac{O(n,n)}{O(n)\times O(n)}$-system} \cite{Ter09}\par

\ss
We use the same notations as in Example \ref{br}, i.e, 
$G=O(n,n,\C)$, and 
$$\tau(g)=\bar g, \quad\sigma_1(g)= \I_{n,n}g\I_{n,n}^{-1}, \quad\sigma_2(g)= \I_{n+1, n-1} g\I_{n+1, n-1}^{-1}.$$ Let $\calA$ be the maximal abelian subalgebra in $\calP_1$ spanned by
$$a_i=\frac{1}{2}\bpm 0 & e_{ii}\\ e_{ii} &0\epm, \quad 1\leq i\leq n,$$
 Then $\calK_1'= 0\times o(n)$, $\calS_1= o(n)\times 0$, and the Lax pair $\o_\l$ of the $\frac{O(n,n)}{O(n)\times O(n)}$-system twisted by $\sigma_2$ is \eqref{dn} with 
 $$g=\bpm \I & 0 \\ 0 & A\epm:\R^n\to K_1', \qquad v_i=\bpm u_i & 0\\ 0 & 0\epm: \R^n\to \calS_1, \quad 1\leq i\leq n.$$
 In other words,
 \beq\label{bg}
 \o_\l =\frac{\l}{2} \bpm 0& \d A^t\\ A\d & 0\epm + \bpm u &0\\ 0&0\epm +\frac{\l^{-1}}{2} \bpm 0& \d A^t J\\ JA\d& 0\epm\,
 \eeq
where $A:\R^n\to O(n)$, $\d=\diag(dx_1, \ldots, dx_n)$, $J=\diag(1, -1, \ldots, -1)$, and $u=\sum_{i=1}^n u_i dx_i$. 

 The flatness of $\o_\l$ is equivalent to $(A, u)$ satisfying the following system
 \beq\label{do}
 \bca 
 \rd A\wedge \d + A\d \wedge u=0,\\
 \rd u + u\wedge u + \d A^t (\frac{\l \I}{2} + \frac{\l^{-1}J}{2})^2 A\d=0.
 \eca.
 \eeq
 The first equation implies that there exists $F=(f_{ij})$ with $f_{ii}=0$ for all $1\leq i\leq n$ such that $$A^{-1}\rd A= \d F^t- F\d, \qquad u=\d F- F^t \d.$$
Since this is the Lax pair \eqref{bg} for the GSGE, the twisted $\frac{O(n,n)}{O(n)\times O(n)}$-system is the GSGE. 
 \eeg
 
\ss\ni{\bf The $-1$ flow on the $\frac{U}{K}$-system}
 
 We combine the $-1$ flow and the flows in the $\frac{U}{K}$-hierarchy generated by $a_i\l$ for $1\leq i\leq n$ to get the $-1$ flow on the $\frac{U}{K}$-system. This is the equation for $v:\R^{n+1}\to \calA^\perp\cap \calP$ and $g:\R^{n+1}\to K$:
 \beq\label{ds}
\bca
 -[a_i, v_{x_j}] + [a_j, v_{x_i}] + [[a_i, v], [a_i, v]] =0, & i\not=j,\\
 [g^{-1}g_{x_i} - [a_i, v], g^{-1}bg]=0, & 1\leq i\leq n,\\
 [a_i, v_t]= [a_i, g^{-1}b g], & 1\leq i\leq n.
 \eca
 \eeq 
  Equation \eqref{ds} has a Lax pair
 $$\o_\l= \left(\sum_{i=1}^n (a_i \l + [a_i, v])\rd x_i\right) + \l^{-1} g^{-1} b g\rd t. $$
 If $\calK_b=\{k\in \calK\n [k, b]=0\}=0$, then the second equation of \eqref{ds} gives $g^{-1}g_{x_i} =[a_i, v]$ for $1\leq i\leq n$. 
 If $\frac{U}{K}$ has maximal rank and $a\in \calA$ is regular, then the $-1$ flow
 on the $\frac{U}{K}$-system becomes the following system:
\beq\label{dsa}
 \bca
 -[a_i, v_{x_j}] + [a_j, v_{x_i}] + [[a_i, v], [a_i, v]] =0, & i\not=j,\\
 g^{-1}g_{x_i}= [a_i, v], & 1\leq i\leq n,\\
[a, v_t]= [a, g^{-1}b g].
 \eca
\eeq
 Note that
 \ben
 \item when $\frac{U}{K}$ is of rank one, the $-1$ flow on the $\frac{U}{K}$-system is the $-1$ flow for the $\frac{U}{K}$-hierarchy by changing the dependent variable $u= [a, v]$,
 \item \eqref{dsa} is an evolution equation on the space of solutions of the $\frac{U}{K}$-system. 
 \een

\ss\ni {\bf Higher flows on the space of solutions of the $\frac{U}{K}$-system}
 
 Assume $a\in \calA$ is a regular element, and $Q=\sum_{j\leq 1} Q_j\l^j$ is constructed from \eqref{bh} using $u=[a,v]$. Note that $Q$ satisfies the recursive formula
 $$(Q_j)_x+ [[a, v], Q_j]= [Q_{j-1}, a],$$
$Q_1=a$, and $f_j(Q)= f_j(a\l)$, where $f_1, \ldots, f_n$ are a set of free generators of the ring of Ad$(K)$-invariant polynomials on $\calP$.  The flow in the $\frac{U}{K}$-hierarchy generated by $a\l^j$ written in $v$ is
\beq\label{en}
[a, v_t]= (Q_{1-j})_x+ [[a, v], Q_{1-j}]= [Q_{j}, a].
\eeq
Recall that $v$ is a solution of the $\frac{U}{K}$-system if and only if $[a, v(x_1, \ldots, x_n)]$ solves the flow generated by $a_i\l$ in the $\frac{U}{K}$-hierarchy for $1\leq i\leq n$.  Since all flows in the $\frac{U}{K}$-system commute, the space of solutions of the $\frac{U}{K}$-system is invariant under the evolution equation \eqref{en} for all odd $j$. In other words,  the following system for $v:\R^n\times \R\to \calA^\perp\cap \calP$,
\beq\label{eo}
\bca -[a_i, v_{x_k}] + [a_k, v_{x_i}] + [[a_i, v], [a_k, v]]=0, & 1\leq i, k\leq n,\\
[a, v_t]= (Q_{1-j})_{x_i} + [[a_i, v], Q_{1-j}], & 1\leq i\leq n,
\eca
\eeq 
has a Lax pair 
$$ (a\l^j + Q_0\l^{j-1} + \cdots + Q_{1-j})\rd t+ \sum_{i=1}^n (a_i\l + [a_i, v])\rd x_i.$$
System \eqref{eo} can be viewed as an evolution equations on the space of solutions of $\frac{U}{K}$-system as follows: Write $a=\sum_{i=1}^n c_i a_i$ and $u=[a,v]$.  Then 
$$\bca 
u_{x_i}=  \ad(a_i)\ad(a)^{-1}(\sum_{i=1}^n c_i u_{x_i} ) +[u, \ad(a_i)\ad(a)^{-1}(u) ], & 1\leq i\leq n,\\
u_t= \sum_{i=1}^n c_i (Q_{1-j}(u))_{x_i}  + [u, Q_{1-j}(u)],
\eca$$
are commuting flows for $u$.  The first set of equation in the above system means that $v=[a,u]$ is a solution of the $\frac{U}{K}$-system.  Hence 
$$[a,v_t]= \sum_{i=1}^n c_i (Q_{1-j}(u))_{x_i}  + [u, Q_{1-j}(u)]$$
leaves the space of solutions of the $\frac{U}{K}$-system invariant.

\bs
\section{Loop group actions} 

 We review the dressing action of $L^{\tau, \sigma}_-(G)$ on the space of solutions of the $\frac{U}{K}$-system, and explain the relation between the action of ``simple'' rational elements in $L^{\tau, \sigma}_-(G)$ and geometric B\"acklund and Ribaucour transforms. 
 
  Let $v$ be a solution of the $\frac{U}{K}$-system, and $E$ the normalized parallel frame for the Lax pair $\o_\l= \sum_{i=1}^n (a_i \l + [a_i, v])\rd x_i$, i.e., $E(x,\l)$ is the solution of 
$$\bca E^{-1}E_{x_i}= a_i\l + [a_i, v],  & 1\leq i\leq n,\\ E(0,\l)=\I. \eca$$
Since $\o_\l$ is holomorphic in $\l\in\C$ and satisfies the $\frac{U}{K}$-reality condition 
$$\tau(\o_{\bar\l}) = \o_\l, \qquad \sigma(\o_{-\l}) = \o_\l,$$
its frame $E(x)\in L^{\tau,\sigma}_+(G)$, where $E(x)(\l)= E(x,\l)$. Given $f\in L_-^{\tau,\sigma}(G)$, by the Local Factorization Theorem \ref{bv}, we can factor 
$$f E(x)= \ti E(x) \ti f(x)$$
with $\ti E(x)\in L_+^{\tau,\sigma}(G)$ and $\ti f(x)\in L_-^{\tau, \sigma}(G)$ in an open subset of $x=0$ in $\R^n$.  Expand 
$$\ti f(x)(\l)= \I + f_1(x)\l^{-1} + \cdots.$$ Then $f_1(x)\in \calP$ and we have

\bthm \label{dp}\cite{TU00a} 

Let $f, v, E, \ti f, \ti f_1, \ti E$ be as above. Then 
\ben
\item $\ti v(x):= (f_1)_\ast$ is a solution of the $\frac{U}{K}$-system, where $(f_1)_\ast$ denotes the projection of $f_1\in \calP$ onto $\calA^\perp\cap \calP$ along $\calA$.
\item $\ti E$ is the normalized parallel frame for $\ti v$.
\item $f\ast v : = \ti v$ defines an action of $L_-^{\tau,\sigma}(G)$ on the space of solutions of the $\frac{U}{K}$-system. 
\item $f\ast E: =\ti E$ defines an action of $L_-^{\tau,\sigma}(G)$ on normalized parallel frames of solutions of the $\frac{U}{K}$-system. 
\item  If $f\in L^{\tau,\sigma}_-(G)$ is rational, then $f\ast v$ can be computed explicitly using $E$ and the poles and residues of $f$. 
\item If $U$ is compact, $a_1$ is regular and $f\in L_-^{\tau,\sigma}(G)$ is rational, then $f\ast 0$ is globally defined and rapidly decaying as $|x_1|\to \infty$.  
\een
\ethm

\brem\hfil\par

\ben
\item
We say $f:S^1\to \C^*\times G$ satisfies the {\it $\frac{U}{K}$-reality condition up to scalar functions\/} if there is a $\phi:S^1\to \C$ such that
$$\tau(f(\bar\l))= \phi(\l) f(\l), \quad \sigma(f(-\l))= \phi(\l) f(\l).$$
Since  scalar functions commute with $L^{\tau, \sigma}(G)$, Theorem \ref{dp} works for rational maps $f$ that satisfy the $\frac{U}{K}$-reality condition up to scalar functions.  
\item Given $f\in L_-^{\tau,\sigma}(G)$, if $E$ is a parallel frame of a solution $v$ of the $\frac{U}{K}$-system and $fE(0, \cdot)$ lies in the big cell of $L^{\tau,\sigma}(G)$ then Theorem \ref{dp}(1) still holds and $\ti E$ is a parallel frame for $f\ast v$ (but may not be normalized).
\een
\erem

\ss\ni {\bf B\"acklund transformations for $U(n)$-system} \cite{TU00a}\par

We use the $U(n)$-system as an example to demonstrate how to compute explicitly the action of the subgroup $\calR_-^{\tau}(G)$ of rational elements in $\calL_-^{\tau}(G)$. Note that $\calR_-^\tau(G)$ is the group of rational maps $f:S^2\to GL(n,\C)$ that satisfying the $U(n)$-reality condition and $f(\infty)=\I$.  
First we find a rational element $f\in \calR_-^\tau(G)$ with only one simple pole, then  use residue calculus to compute the action of $f$ on solutions of the $U(n)$-system.   

Let  $\a\in \C$, $\pi$ a Hermitian projection of $\C^n$, and $\pi^\perp=\I-\pi$. Then
\beq\label{cu}
g_{\a, \pi}(\l) = \pi + \frac{\l-\bar\a}{\l -\a} \pi^\perp = \I +\frac{\a-\bar\a}{\l-\a} \pi^\perp
\eeq
satisfies the $U(n)$-reality condition $g(\bar\l)^*g(\l)=\I$.

\ss
\ni {\bf Three methods to compute $g_{\a, \pi}\ast v$}\par

\ss \ni {\it Method 1: Algebraic B\"acklund Transformation}\par 

\ss
Let $\calA$ be the space of diagonal matrices in $u(n)$, $a_j=\im e_{jj}$,  $v$ a solution of the $U(n)$-system, and $E$ the normalized parallel frame, i.e., $E^{-1}\rd E= \o_\l=\sum_{i=1}^n (a_i\l + [a_i, v])\rd x_i$ and $E(0,\l)=\I$. We claim that
$$g_{\a, \pi}\ast v= v+ (\a-\bar\a)\ti\pi_\ast,$$ where $\ti\pi(x)$ is the Hermitian projection of $\C^n$ onto $E(x, \a)^{-1}(\Im \pi)$ and $\ti\pi_\ast$ is the projection of $u(n)$ onto $\calA^\perp$ along $\calA$. To see this, we need to factor $g_{\a,\pi} E(x) = \ti E(x) \ti g(x)$ with $\ti E(x)\in L^\tau_+(G)$ and $\ti g(x)\in L^\tau_-(G)$. We make an Ansatz that $\ti g= g_{\a, \ti\pi(x)}$ and  solve $\ti\pi(x)$ by  requiring that
\begin{align*}
\ti E(x,\l): &= g_{\a, \pi}(\l)E(x,\l) \ti g^{-1}(x,\l)\\
&= (\I+ \frac{\a-\bar\a}{\l- \a}\pi^\perp) E(x,\l) (\I -\frac{\a-\bar\a}{\l- \bar\a}\ti\pi(x)^\perp)
\end{align*}
lies in $L^\tau_+(G)$.  Hence the residues of $\ti E(x,\l)$ at $\l=\a, \bar\a$ should be zero. This implies that
$$\pi^\perp E(x, \a) \ti \pi(x)=0, \qquad \pi E(x, \bar\a) \ti \pi(x)^\perp=0.$$
Both conditions are satisfied if $$\Im (\ti\pi(x))= E(x,\a)^{-1}(\Im (\pi)).$$  This gives the formula for $\ti\pi (x)$.  The formula for $\ti E$ implies that $\ti E^{-1}\rd \ti E$ has a simple pole at $\l=\infty$ and $\ti E^{-1}\rd \ti E= \sum_{i=1}^n a_i\l + [a_i, \ti v]$, where $\ti v= v+ (\a-\bar\a)\ti\pi_\ast$. This proves the claim.

\ms\ni {\it Method 2: ODE B\"acklund transformation}\par

The new solution $g_{\im s, \pi}\ast v$ can be also obtained by solving a system of compatible ODEs: 
Set $\o_\l= E^{-1} \rd E$ and $\ti \o_\l= \ti E^{-1}\rd \ti E$.
Since $\ti E= g_{\a,\pi} E g_{\a, \ti \pi}^{-1}$ and $g_{\a, \pi}$ is independent of $x$, $\ti \o_\l = \ti g \o_\l \ti g^{-1} - \rd \ti g \ti g^{-1}$;  or equivalently,
\beq\label{cr}
\ti \o_\l \ti g= \ti g\o_\l - \rd \ti g,
\eeq
 where $\ti g= g_{\a, \ti\pi}$.
Multiply  \eqref{cr} by $(\l-\a)$ and compare coefficients of $\l^i$ to see that $\ti \pi$ must satisfy
\beq \label{ce}
\bca \ti\pi_{x_j}+  [\a a_j+ [a_j, v], \, \ti \pi] = (\a-\bar\a) [a_j, \ti\pi]\ti\pi^\perp, & 1\leq j\leq n,\\ \ti\pi(x)^*= \ti \pi(x).\quad \ti \pi^2=\ti \pi,\eca 
  \eeq  
  and $g_{\a, \pi}\ast v= v+ (\a-\bar\a)\ti\pi_\ast$
Moreover, given $v$,
\ben 
\item system \eqref{ce} is solvable for $\ti\pi$ if and only if $v$ is a solution of the $U(n)$-system,
\item if $v$ is a solution of the $U(n)$-system and $\ti\pi$ the solution of \eqref{ce}, then  $\ti v= v+ (\a-\bar\a) \ti \pi_\ast$ is a solution of the $U(n)$-system, where $\ti\pi_\ast$ is the projection of $\ti\pi$ onto $\calA^\perp$ along $\calA$.  
 \een
 
\ms\ni{\it Method 3: Linear B\"acklund transformations}\par

Suppose $\pi$ is the Hermitian projection of $\C^n$ onto $V= \C y_0$. Set 
$$y(x)= E(x, \a)^{-1}(y_0).$$  The normalized parallel frame of $g_{\a, \pi}\ast v$ is 
$g_{\a, \pi} E(x,\cdot) g_{\a, \ti\pi(x)}^{-1}$, where $\ti\pi(x)$ is the projection onto $\C y(x)$.   Differentiate $y$ to get 
$$\rd y= -E^{-1}\rd E E^{-1} y_0= -\o_\a y.$$  So 
$y$ is the solution of the following linear system
\beq\label{bj}
\bca \rd y= -\o_\a y= -\sum_{j=1}^n ( a_j\a +[a_j, v]) \rd x_j,\\
  y(0)= y_0.\eca
\eeq
In fact, given $v:\R^n\to \calA^\perp\cap\calP$, 
\ben
\item system \eqref{bj} is solvable if and only if $v$ is a solution of the $\frac{U}{K}$-system,
\item if $v$ is a solution of the $\frac{U}{K}$-system and $y$ is a solution of \eqref{bj}, then $g_{\a,\pi}\ast v= v+ (\a-\bar\a)\ti\pi_\ast$, where $\ti\pi(x)$ is the Hermitian projection of $\C^n$ onto $\C y(x)$.  
\een

\ss
Note that the first and third methods are essentially the same because solutions $y$ of \eqref{bj} is $E(\cdot, \a)^{-1}(y_0)$, where $E(\cdot, \a)$ is a parallel frame for $\o_\a$. 

If $\dim(\Im \pi)=k$, then we first choose a basis $y^0_1, \ldots, y^0_k$ of $\Im\pi$. Let $y_i$ be the solution of  \eqref{bj} with $y_i(0)= y_i^0$, $\ti V(x)$ the linear subspace spanned by $y_1(x), \ldots, y_k(x)$, and $\ti \pi(x)$  the Hermitian projection  of $\C^n$ onto $\ti V(x)$.  Then the new solution is $g_{\a,\pi}\ast v= v+ (\a-\bar\a)\ti\pi_\ast$.  

\ss\ni {\bf Permutability formula} \cite{TU00a}

The permutability formulae for B\"acklund transformations for the SGE, the GSGE, Ribaucour transforms for flat Lagrangian submanifolds in $\C^n$ and  for $k$-tuples in $\R^n$ of type $\R^{k-\ell, \ell}$ can be obtained in a unified way. This is because 
\ben
\item geometric transforms on these submanifolds correspond to actions of simple rational elements in the negative loop group, 
\item  if $g_i$ have poles at $\a_i$ for $i=1,2$, then we use residue calculus to factor $g_1g_2= f_2f_1$ such that $f_i$ have poles at $\a_i$ for $i=1,2$. 
\een 
Permutability formulae can then be obtained from the fact that the geometric transforms are actions. 

   We use $U(n)$-system as an example to explain this method:
Given  $g_{\a_1,\pi_1}$, $ g_{\a_2, \pi_2}$ with $\a_1\not= \pm\bar\a_2$, let $\tau_1, \tau_2$ be the projections  such that 
\beq\label{ch}
\Im \tau_1= g_{\a_2,\pi_2}(\a_1)(\Im \pi_1), \quad \Im\tau_2= g_{\a_1, \pi_1}(\a_2)(\Im \pi_2).
\eeq
Then
\beq\label{cg}
g_{\a_2, \tau_2} \circ g_{\a_1, \pi_1} = g_{\a_1, \tau_1} \circ g_{\a_2, \pi_2}.
\eeq
This gives a relation for rational elements in $R_-^\tau(G)$ with only one simple pole. 

Formula \eqref{cg} leads to a Bianchi type permutability formulae for B\"acklund transformations as follows:  
Let $v_0$ be a solution of the $U(n)$-system, and $E_0(x,\l)$ its normalized parallel frame.  Let $\ti \pi_j(x)$ denote the Hermitian projections of $\C^n$ onto $E_0(x,\a_j)^{-1}(\Im\pi_j)$ for $j=1,2$.  Then 
$$E_j(x,\l) = g_{\a_j, \pi_j}(\l) E_0(x,\l) g_{\a_j, \ti\pi_j(x)}(\l)^{-1}$$
is the normalized parallel frame for 
$$v_j= g_{\a_j, \pi_j}\ast v_0= v_0 + (\a_j-\bar\a_j) (\ti\pi_j)_\ast, \quad j= 1, 2. $$
Use the fact that $L_-^\tau(G)$ acts on the space of solutions and the permutability formula \eqref{cg} to get
$$v_3=g_{\a_2,\tau_2}\ast v_1= g_{\a_2, \tau_2}\ast (g_{\a_1, \pi_1}\ast v_0)= g_{\a_1, \tau_1} \ast(g_{\a_2, \pi_2}\ast v_0)= g_{\a_1, \tau_1}\ast v_2.$$
But 
\begin{align*}
&v_3= v_1 + (\a_2-\bar\a_2) (\ti \tau_2)_\ast =v_2 + (\a_1-\bar\a_1) (\ti \tau_1)_\ast, \quad {\rm where\/}\\
&\Im\ti\tau_2=E_1(x, \a_2)^{-1}(\Im \tau_2), \quad \Im\ti\tau_1= E_2(x, \a_1)^{-1}(\Im \tau_1).
\end{align*}
So $v_3$ can be given by an explicit formula in terms of $v_0, v_1, v_2$.  This gives the {\it permutability formula} for the $U(n)$-system.

\ss\ni {\bf Action of $\calR_-^\tau(G)$}

\ss
The method we used to construct the action of $g_{\a, \pi}\ast v$ works for the action of any $f\in \calR_-^\tau(G)$ on $v$ as follows:  First we write 
$$f(\l)= \I + \sum_{i=1, j=1}^{k, n_i} \frac{P_{ij}}{(\l-\a_i)^j}$$
for some constants $\a_i\in \C$ and $P_{ij}\in gl(n)$.  Let $E$ be the normalized parallel frame of a solution $v$ of the $\frac{U}{K}$-system.  We assume $fE(x)= \ti E(x) \ti f(x)$ where $\ti f(x)$ has poles at $\a_1, \ldots, \a_k$ with order $n_1, \ldots, n_k$ respectively, i.e., 
$$\ti f(x,\l)= \I + \sum_{i=1, j=1}^{k, n_i} \frac{\ti P_{ij}(x)}{(\l-\a_i)^j}$$
  Reality condition gives $\ti f(x,\l)^{-1}= \overline{\ti f(x,\bar\l)}^t$. Then $f(\l) E(x, \l) \ti f(x, \l)^{-1}= f(\l) E(x,\l) \ti f(x, \bar\l)^*$ should have no poles at $\l=\a_i$ for $1\leq i\leq k$. We can use these conditions to solve $\ti P_{ij}(x)$.  This computation is long and tedious.  However, if we find a set of generators of the negative rational loop group $\calR_-^\tau(G)$ with minimal number of poles then we can simplify the computation by using  permutability formulas (relations) for these generators or the algebraic BT.

\ss\ni {\bf Simple elements and generators}

Let $\frac{U}{K}$ denote the symmetric space constructed from two commuting involutions $\tau, \sigma$, and $\calR_-^{\tau,\sigma}(G)$ denote the subgroup of rational maps $f:S^2\to G$ that are in  $ L_-^{\tau,\sigma}(G)$.  
  A $f\in \calR_-^{\tau, \tau}(G)$ is called {\it a simple element\/} if $f$ can not be factored as product of $f_1f_2$ with both $f_1$ and $f_2$ in $\calR_-^{\tau,\sigma}(G)$.   The following are known:

\ben
\item Uhlenbeck \cite{Uhl89} proved that 
$$\{g_{\a, \pi}\n \a\in \C, \pi^*=\pi, \pi^2=\pi\}$$ generates the negative rational loop group satisfying the $U(n)$-reality condition.

\item Note that
\ben
\item $g_{\im s, \pi}$ satisfies the $\frac{U(n)}{O(n)}$ reality condition if  $s\in \R$ and $\bar\pi=\pi$.
\item if $\a\in \C\setminus \im \R$, $\pi$ is a Hermitian projection of $\C^n$, and $ \Im\rho= g_{\a,\pi}(-\bar \a)(\Im\bar\pi)$, then $$f_{\a,\pi}=g_{-\bar\a, \rho}g_{\a, \pi}$$ satisfies the $\frac{U(n)}{O(n)}$ reality condition. 
\een
 Terng and Wang \cite{TW08} proved that these elements generate the negative rational loop group satisfying the $\frac{U(n)}{O(n)}$-reality condition.
\item Donaldson, Fox, and Goertsches \cite{DFG08} construct a set of generators for $\calR^{\tau,\sigma}(G)$ when $G$ is a classical group.
\een

\ss\ni{\bf B\"acklund transforms for $\frac{U(n)}{O(n)}$-system} \cite{TU00a}

The methods described above for constructing algebraic and analytic BT and permutability formula for $U(n)$-system work the same way for general $\frac{U}{K}$-system. For example,  $g_{\im s, \pi}$ satisfies the $\frac{U(n)}{O(n)}$ reality condition. If
 $v$ is a solution of the $\frac{U(n)}{O(n)}$-system and $E$ is its normalized parallel frame for the Lax pair of $v$, then:
\ben
\item $E(x, \cdot)$ satisfies the $\frac{U(n)}{O(n)}$ reality condition.
\item Since $\o_{\im s}= \sum_{i=1}^n (\im s a_i + [a_i, v])\rd x_i$ and $a_i$ is diagonal in $u(n)$, $\o_{\im s}$ is a $sl(n,\R)$-valued $1$-form.  Hence $E(x,\im s)\in SL(n, \R)$ and $E(x, \im s)(\Im \pi)$ is real.
\item $g_{\im s, \pi}\ast v= v+ 2\im s \ti \pi_\ast$ is a solution of the $\frac{U(n)}{O(n)}$-system, where $\ti \pi$ is the orthogonal projection of $\R^n$ onto $E(x, \im s)^{-1}(\Im \pi)$.  
\een

\bs
\section{Action of simple elements and geometric transforms}

Suppose a class of submanifolds in Euclidean space admits a local coordinate system and an adapted frame such that its Gauss-Codazzi equation is the $\frac{U}{K}$-system (or twisted $\frac{U}{K}$-system) for some symmetric space $\frac{U}{K}$.  If the adapted frame and the immersion of the submanifold can be obtained from the parallel frame of the Lax pair of the corresponding solution of the $\frac{U}{K}$-system, then the action of a simple rational loop on the parallel frame of a solution of the $\frac{U}{K}$-system gives rise to a geometric transform of these submanifolds.  We explain how this is done for $K=-1$ surfaces in $\R^3$, flat Lagrangian submanifolds in $\C^n$, and $k$-tuples in $\R^n$ of type $\R^{k-\ell, \ell}$.  We have given a unified method to construct Permutability formula for actions of simple elements on the space of solutions and normalized parallel frames of $\frac{U}{K}$-systems in section 9. Hence if we know how to read geometric transforms from the action of simple elements on parallel frames then we can obtain an analogue of Bianchi's Permutability Theorem for these geometric transforms.   

\ss\ni {\bf BT for $K=-1$ surfaces in $\R^3$ and action of $g_{\im s, \pi}$} 

Let $g_{\im s,\pi}$ be the rational map defined by \eqref{cu} with $s\in \R$ and $\pi$ real.
It was noted by Uhlenbeck in \cite{Uhl92} that the dressing action of $g_{\im s,\pi}$ on solutions SGE gives rise the B\"acklund transforms for $K=-1$ surfaces in $\R^3$. 

Let $q$ be a solution of the SGE, $2q_{xt}= \sin 2q$, and   $E(x,t,\l)$ the normalized parallel frame for the Lax pair 
$$\o_\l= \left(\l \bpm -\im &0\\ 0& \im\epm + \bpm 0 & -q_x\\ q_x &0\epm \right) \rd x + \frac{i}{4\l} \bpm \cos 2q & -\sin 2q\\ -\sin 2q & -\cos 2q\epm \rd t.$$
Then 
\beq\label{dz}
f= \frac{\p E}{\p \l}E^{-1}\bigg|\, _{\l= \frac{1}{2}}
\eeq
 is the immersion of a $K=-1$ surface in $su(2)$ (identified as $\R^3$) corresponding to the solution $q$ of SGE.  We have seen that $\ti E= g_{\im s, \pi} E g_{\im s, \ti \pi}^{-1}$ is the normalized parallel frame for $g_{\im s, \pi}\ast q$, where $\ti \pi(x)$ is the orthogonal projection of $\R^2$ onto $E(x,\im s)^{-1}(\Im \pi)$.  Then 
\beq\label{dx}
\hat E=\left(\frac{\l+\im s}{\l -\im s}\right)^{\frac{1}{2}}Eg_{\im s, \ti \pi}^{-1}.
\eeq
is a parallel frame for $g_{\im s, \pi}\ast q$, and 
\beq\label{dy}
\hat f= \frac{\p \hat E}{\p \l}\hat E^{-1}\,\bigg|_{\l= \frac{1}{2}}
\eeq
is the immersion of a $K=-1$ surface in $su(2)$ corresponding to $g_{\im s, \pi}\ast q$.  Note that $\hat E\in SU(2)$.   To see the properties of the transform $f\mapsto \hat f$, we use \eqref{dx} and \eqref{dy} to get 
$$
\hat f= f + \frac{2\im s}{\frac{1}{4} + s^2}E(\cdot,\frac{1}{2})(\ti\pi^\perp -\frac{1}{2}\I) E(\cdot,{\frac{1}{2}})^{-1}.
$$
Let $(\cos y(x), \sin y(x))^t$ denote the unit direction of the real line $\Im\ti\pi(x) \subset \R^2$.  Then a direct computation then implies that
$$\hat f= f+ \sin \o e_1,$$ where $\sin\o= \frac{s}{\frac{1}{4} + s^2}$ and 
$$e_1= \cos 2y\, E_{\frac{1}{2}}\bpm -\im & 0\\ 0 & \im\epm E_{\frac{1}{2}}^{-1} - \sin 2y \, E_{\frac{1}{2}} \bpm 0& \im \\ \im & 0\epm E_{\frac{1}{2}}^{-1}$$
is tangent to $f$, where $E_{\frac{1}{2}}=E(\cdot, \frac{1}{2})$. 
Use  \eqref{dx} to see that $\hat f-f$ is tangent to $\ti f$. In other words, $f\mapsto \hat f$ is a BT with angle $\o$.

\ss\ni {\bf $n$-submanifolds in $\R^{2n-1}$ with constant curvature $-1$}

Let $L_\pm$ denote the positive and negative groups defined in Example \ref{ct} for the $\frac{O(n,n)}{O(n)\times O(n)}$-system twisted by $\sigma_2$. First we construct a simple rational map satisfies the  $\onn$-reality condition up to scalar functions.  A direct computation shows that if $g(\l)= \bpm 1& 0\\ 0& \b\epm +\frac{s}{\l- s} P$ satisfies the $\frac{O(n,n)}{O(n)\times O(n)}$-reality condition up to scalar functions, then 
$$P=\bpm 1& C^t\\ C& 1\epm \bpm 1&0\\ 0& \b\epm.$$
In other words, $g$ must be of the form
$$g_{\b, C}(\l)= \bpm 1&0\\ 0& \b\epm + \frac{s}{\l-s} \bpm 1& C^t\b\\ C& \b\epm = \frac{1}{\l-s} \bpm\l & sC^t\b\\ sC & \l \b\epm,$$
where $\b, C\in O(n)$.

Let $A$ be a solution of the GSGE, and $E(x,\l)$ the normalized parallel frame for the corresponding Lax pair $\o_\l$ defined by \eqref{bs}.  Note that $E(x,\cdot)\in L_+$.
Suppose $g_{\b,\I}(\l)E(x,\l)= \ti E(x,\l) g_{\ti\b(x), \ti C(x)}(\l)$
with $\ti\b(x), C(x)\in O(n)$ and $\ti E(x,\cdot)\in L_+$.  Then 
$$\ti E(x,\l)= g_{\b, \I}(\l) E(x, \l ) g_{\ti \b(x), \ti C(x)}(\l)^{-1}$$ is holomorphic for $\l\in \C$.  So the residue at $\l=s$ is zero, i.e.,
$$(1,\b)E(x, s) \bpm 1 & -\ti C^t\\ -\ti \b^t \ti C & \ti \b^t \epm =0.$$
This implies that 
$$\ti\b^t\ti C= (\eta_2 + \b \eta_4)^{-1} (\eta_1 + \b \eta_3), \quad {\rm where\,\,} E(\cdot, s)=\bpm \eta_1& \eta_2\\ \eta_3 & \eta_4\epm.$$
Set
$$(P,Q):= (1, \b) E(\cdot, s).$$
Then 
$\rd (P, Q) = (P, Q) \o_s= (P, Q) \bpm w & \d A^t D_s\\ D_sA \d & 0\epm, \quad D_s= \frac{1}{2}(s\I + s^{-1}J)$, or equivalently,
$$\bca \rd P = P w+ QD_s A\d,\\ \rd Q= P\d A^t D_s.\eca$$
If $ X:= -Q^{-1}P$, then we get the BT given in Theorem \ref{bd}:
$$\rd X= X\d A^t D_s X - Xw - D_s A\d.$$
This explains the following Theorem of \cite{BT88} in terms of the action of $g_{\b, C}$:

\bthm
Let $s$ be a non-zero real constant. Consider the linear system for $y:\R^n\to \calM_{n\times 2n}$:
\beq\label{ff}
\rd y= y \bpm w & \d A^t D_s\\ D_sA \d & 0\epm, \quad D_s= \frac{1}{2}(s\I + s^{-1}J).
\eeq
Then
\ben
\item System \eqref{ff} is solvable if and only if $A$ is a solution of the GSGE.
\item   If $y=(P, Q)$ is a solution of \eqref{ff} with $Q\in GL(n)$,  then 
 $X= -Q^{-1}P$ is a solution of BT \eqref{em} for GSGE and $X$ is a solution of GSGE.  
 \een
\ethm

In other words, \eqref{ff} can be viewed as the {\it Linear B\"acklund transform} for GSGE.

\bdefn  {\bf Ribaucour transform for submanifolds}  \cite{DT02}\par\hfil

Let $M$ and $\ti M$ be two $n$-dimensional submanifolds in $\R^{n+k}$ with flat normal bundle. A {\it Ribaucour transform\/} is a vector bundle isomorphism $\Phi:\nu(M)\to \nu(\ti M)$ covers a diffeomorphism $\phi:M\to \ti M$ satisfying the following conditions:
\ben
\item $\Phi$ maps parallel normal fields of $M$ to parallel normal fields of $\ti M$,
\item for each $p\in M$ and $v\in \nu(M)_p$, the normal line $p+ tv$ intersects the normal line $\phi(p)+ t \Phi(v)$ at equal distance $r(p,v)$,
\item $\rd\phi_p$ maps common eigenvectors of shape operators of $M$ at $p$ to common eigenvectors of shape operators of $\ti M$ at $\phi(p)$,
\item the tangent line through $p$ in a principal direction $v$ meets the tangent line through $\phi(p)$ in the direction of $\rd\phi_p(v)$ at equal distance,
\een
\edefn

Let $M$ be a submanifold in $\R^{n+k}$, and $(e_1, \ldots, e_{n+k})$ an orthonormal frame on $M$ such that $(e_1, \ldots, e_n)$ are principal directions (i.e., a common eigen-frame for the shape operator of $M$) and $(e_{n+1}, \ldots, e_{n+k})$ is a parallel normal frame.  Let $\ti M$ be another $n$-submanifold with flat normal bundle, $\phi:M\to \ti M$  a diffeomorphism, $(\ti e_{n+1}, \ldots, \ti e_{n+k})$ a parallel normal frame for $\ti M$, and $\ti e_i$ is the direction of $\rd \phi(e_i)$ for $1\leq i\leq n$.    Then $\phi$ is a Ribaucour transform if 
\ben
\item[(a)] $\ti e_i$ is a principal direction for $\ti M$ for $1\leq i\leq n$,
\item[(b)]
 there exist functions $h_1, \ldots, h_{n+k}$ on $M$ such that
\[\phi(p)+ h_i(p)\ti e_i (p)= p+ h_i (p) e_i(p), \qquad 1\leq i\leq n+k\]
for all $p\in M$.   
 \een  

\ss\ni {\bf Flat Lagrangian submanifolds in $\C^n$} 

Let $(\b, h)$ be a solution of the $\frac{U(n)\ltimes \C^n}{O(n)\ltimes \R^n}$-system, $\o_\l$ its Lax pair \eqref{dm},  and $F=\bpm E & X\\ 0&1\epm$ the normalized parallel frame of $\o_\l$.   We have seen in section \ref{da} that for each $r\in \R$, $X(\cdot, r)$ is a flat Lagrangian immersion in $\C^n$ corresponding to solution $(\b, h)$ (the associated family). We review the action of two types of simple elements on the space of solutions of the
$\frac{U(n)\ltimes \C^n}{O(n)\ltimes \R^n}$-system and derive the corresponding geometric transformations (\cite{TW08}). 

\ss\ni {\it The action of $h_{\a, \pi}$} 

We compute the action of $h_{\a, \pi}$ on flat Lagrangian submanifolds in $\C^n$, where 
$$h_{\a,\pi} =\bpm g_{\im \a, \pi} &0\\ 0 & \frac{\l+\im \a}{\l -\im \a}\epm$$ with $\a\in\R$ and $\bar\pi=\pi$.  Note that $h_{\a,\pi}$ satisfies the $\frac{U(n)\ltimes \C^n}{O(n)\ltimes \R^n}$ reality condition up to scalar functions. 

We claim that the action of $h_{\a,\pi}$ gives a Ribaucour transform for flat Lagrangian submanifolds in $\C^n$.
To see this, first we factor $gF= \ti F \ti f$ with 
$$\ti F= \bpm \ti E & \ti X\\ 0&1\epm, \qquad \ti f= \bpm g_{\im \a, \ti\pi} & \xi\\ 0 & \frac{\l+\im \a}{\l - \im \a}\epm,$$
where $\xi= \frac{-2\im\a}{\l-\im \a} \ti \pi \eta$, $\eta(x)= E(x, -\im \a)^{-1} X(x, -\im \a)$, and $\ti\pi(x)$ is the Hermitian projection onto $\ti y(x)=E(x, \im \a)^{-1}(\Im \pi)$.  It follows from reality conditions that both $\ti\pi$ and $\eta$ are real.  

 We assume $\Im\pi$ is of one dimension and is equal to $\R y_0$.  Let 
$$\ti y(x) = E(x, \im \a)^{-1}(y_0).$$
Then 
\beq\label{er}
\ti \pi = \frac{\ti y\ti y^t}{||\ti y ||^2}.
\eeq
 Equate the $12$-entry of $fF=\ti F \ti f$ to get
$gX= \ti E\xi + \frac{\l+\im \a}{\l- \im \a} \ti X$. This implies that 
$$ X= \ti g^{-1}\ti E \xi + \frac{\l+\im \a}{\l-\im \a} g^{-1}\ti X,$$
where $g= g_{\im \a, \pi}$ and $\ti g= g_{\im \a, \ti \pi}$.
Set $\hat X=\frac{\l+\im \a}{\l -\im \a} g_{\im\a, \pi}^{-1} \ti X$. Then
\beq\label{ep}
\hat X= X + \frac{2\im \a}{\l-\im \a} E\ti\pi \eta
\eeq
is a flat Lagrangian submanifold in $\C^n$ corresponding to the solution $(\ti \b, \ti h)$, where
$\ti \b= \b -2\a (\ti\pi)_\ast$ and $ \ti h= h - 2\a \ti\pi \eta$.

Claim that \eqref{ep} is a Ribaucour transform.  To see this we first note that 
\beq\label{et}
\hat E= \frac{\l+\im \a}{\l-\im \a} g_{\im\a, \pi}^{-1} \ti E = E(\I + \frac{2\im \a}{\l-\im \a} \ti \pi)
\eeq
is a parallel frame for the Lax pair of $(\ti \b, \ti h)$.  Hence 
\beq\label{es}
\hat E- E= \frac{2\im \a}{\l-\im \a} E\ti \pi.
\eeq
Write $E=(e_1, \ldots, e_n)$ and $\hat E= (\hat e_1, \ldots, \hat e_n)$. 
By \eqref{er}, we see that the $j$-th column of $\frac{2\im \a}{\l-\im \a} E\ti \pi$ is equal to $\ti y_j Z$, where 
$$Z= \frac{2\im \a}{\l-\im\a} \frac{E\ti y}{||\ti y ||^2}.$$ 
By \eqref{ep} and \eqref{es}, we get
\begin{gather}
\hat e_j- e_j = \ti y_j Z, \label{ev}\\
\hat X-X= \frac{(\ti y, \eta)}{\ti y_i} (\hat e_i -e_i) \label{ew}
\end{gather}
It remains to compute the relation between parallel normal fields of $X$ and $\hat X$.  The parallel tangent frames for $X$ and $\hat X$ are $V=EA^{-1}$ and $\hat V=\hat E \hat A^{-1}$ respectively, where $A(x)= E(x,0)$ and $\hat A(x)= \hat E(x,0)$.  By \eqref{et}, $\hat A= A(\I- 2\ti\pi)$.  Compute directly to see that
\begin{align*}
\hat V&= \hat E\hat A^{-1}= E(\I+ \frac{2\im \a}{\l-\im \a} \ti \pi) (\I-2\ti\pi) A^{-1} \\
&= E(\I -\frac{2\l}{\l-\im\a}\ti \pi) A^{-1}= EA^{-1} -\frac{2\l}{\l-\im\a} E\ti \pi A^{-1}\\
&= V- \frac{2\l}{\l -\im\a} E\ti \pi A^{-1}.
\end{align*}
Thus we have
\beq\label{eu}
\hat v_j - v_j= \frac{\im\l}{\a} (\sum_{k=1}^n a_{jk} \ti y_k) Z, 
\eeq
where $\ti v_j$ and $v_j$ are the $j$-th column of $\hat V$ and $V$ respectively. Since $X, \hat X$ are Lagrangian, $v_{n+j}=\im v_j$ and $\hat v_{n+j} = \im \hat v_j$ are parallel normal fields for $X$ and $\hat X$ respectively. As a consequence of \eqref{eu}, \eqref{ev} and \eqref{ew}, we have
$$\hat X- X=  -\frac{\a(\eta, \ti y)}{\l \sum_{k=1}^n a_{jk} \ti y_k} (\hat v_{n+j}- v_{n+j}).$$ 
This proves that $X\mapsto \hat X$ is a Ribaucour transform.  In fact, this is the Ribaucour transform found in \cite{DT00}. 

\ss\ni{\it The action of $k_{\im \a, b}$}

We claim that the action of $k_{\im \a, b}$ gives an Combescure O-transform for flat Lagrangian submanifolds in $\C^n$, where 
$$k_{\im \a, b}(\l)= \bpm \I & \frac{\im b}{\l-\im\a} \\0 & 1\epm.$$  

First factor 
$$k_{\im \a, b}F= \ti F \ti k= \bpm E & Y\\ 0 & 1\epm \bpm \I & \frac{\im E_{\im \a}^{-1}b}{\l-\im \a}\\ 0 & 1\epm. $$ 
Then 
$$Y= X+ \frac{\im (b-E_\l E_{\im \a}^{-1} b)}{\l-\im\a}.$$
Moreover, if $\l\in \R$ then $Y$ is a flat Lagrangian submanifold of $\C^n$ corresponding to the solution $(\b, \ti h)$, where 
$\ti h= h + E(\cdot, \im \a)^{-1}b$.  
Note that the transform $X\mapsto Y$ is a Combescure O-transform.

\ss\ni {\bf $k$-tuples in $\R^n$ of type $\R^{k-\ell, \ell}$}

It is known that the Darboux (or Ribaucour) transforms for Christoffel pairs of isothermic surfaces in $\R^3$ and for Christoffel pairs of isothermic surfaces in $\R^n$ can be derived from the action of a simple rational map by dressing actions  (cf. \cite{C97, HP97} and \cite{BDPT02, Bu04} respectively).  Ribaucour transforms are constructed for $k$-tuples in $\R^n$ of type $\R^{k-\ell, \ell}$ in \cite{BDPT02, DonTer08a} using dressing action of a simple rational loop.  Recall that 
Christoffel pairs of isothermic surfaces in $\R^n$ (for $n\geq 3$) are $2$-tuples in $\R^n$ of type $\R^{1,1}$.  So the construction of Ribaucour transforms for $k$-tuples in $\R^n$ of type $\R^{k-\ell, \ell}$ contains the surface case.

\ss\ni {\it Simple elements  for the $\frac{O(n+k-\ell,\ell)}{O(n)\times O(k-\ell,\ell)}$-system} 

Let $W\in \R^n$ and $Z\in \R^{k-\ell, \ell}$ with length 1, i.e., $W^tW=Z^t\I_{k-\ell, \ell}Z=1$, and $\pi$  the projection of $\C^{n+k}$ onto $\C\bpm W \\ \im Z\epm$, i.e., 
$$\pi=\frac{1}{2} \bpm WW^t& \im WZ^t\\ \im ZW^t & ZZ^t\epm.$$  Note that $\pi\bar\pi = \bar\pi \pi=0$.  Let $s\in \R$. Then 
\begin{align*}
p_{\im s, \pi}&= (\pi +\frac{\l+ \im s}{\l -\im s} (\I-\pi)) (\bar \pi + \frac{\l-\im s}{\l+ \im s} (\I-\bar\pi))\\
&= \frac{\l+ \im s}{\l-\im s}\bar\pi + \frac{\l-\im s}{\l+ \im s}\pi + \I-\pi-\bar\pi
\end{align*}
satisfies the $\frac{O(n+k-\ell, \ell)}{O(n)\times O(k-\ell, \ell)}$-reality condition.   

\bthm\label{eh}  \cite{BDPT02, DonTer08a}\hfil\par 

Let $\xi=(F,\g)$ be a solution of the $\frac{O(n+k-\ell, \ell)}{O(n)\times O(k-\ell,\ell)}$-system, and $E(x,\l)$ a parallel frame for the Lax pair $\o_\l$ defined by \eqref{ef}. Let $W\in \R^{n}$ and $Z\in \R^{k-\ell,\ell}$ be unit vectors, $\pi$ the projection of $\C^{n+k}$ onto $\C\bpm W\\ \im Z\epm$.  Then:
\ben
\item  $E(x, \im s)^{-1}\bpm W\\ \im Z\epm$ is of the form $\bpm \ti W(x)\\ \im \ti Z(x)\epm$ with $\ti W\in \R^{n}$, $\ti Z\in \R^{k-\ell, \ell}$, and $\ti W^t\ti W=\ti Z^t\I_{k-1,1}\ti Z$. 
\item The action $p_{\im s, \pi}\ast \xi= (F, \g)+ 4s (\hat Z\hat W^t)_\ast$, where  $\eta_\ast=\eta-\sum_{i=1}^k \eta_{ii} e_{ii}$ for $k\times n$ matrix $\eta= (\eta_{ij})$ and  $\hat W(x)$ and $ \hat Z(x)$ are the unit directions of $\ti W(x)$ in $\R^{n}$ and $\ti Z(x)$ in $\R^{k-\ell, \ell}$ respectively.
\item $\hat E(x,\l):=E(x, \l) p_{\im s, \hat \pi(x)}(\l)$ is a parallel frame for the Lax pair of $p_{\im s, \pi}\ast \xi$, where $\hat \pi$ is the projection onto $\C(\hat W, \im\hat Z)^t$.
\een 
\ethm

We use Theorem \ref{eh}, $\hat E= Ep_{\im s,\hat\pi}$ and a  straight-forward computation to write down the geometric transform on $k$-tuples of type $\R^{k-\ell, \ell}$ corresponding to the action of $p_{\im s, \pi}$. We state the results for the case $n=k+1$, and similar results hold for higher co-dimension. 

\bthm\label{el}{\bf Ribaucour transform for $k$-tuples} \cite{BDPT02, DonTer08a}\hfil\par

Let $E, p_{\im s, \hat\pi}, \hat E$ be as in Theorem \ref{eh}, and $n=k+1$. Then:
\ben 
\item There are $\calM_{(k+1)\times k}$ valued maps $\Xi, \hat \Xi$ such that
$$\frac{\p E}{\p \l}E^{-1}\, \bigg|_{\l=0}= \bpm 0 & \Xi\\ -\Xi^t J& 0\epm,\qquad \frac{\p \hat E}{\p \l}\hat E^{-1}\, \bigg|_{\l=0}= \bpm 0 & \hat\Xi\\ -J\hat \Xi^t & 0\epm,$$ 
where $J=\I_{k-\ell,\ell}=\diag(\e_1, \ldots, \e_k)$.
\item Given a non-zero vector $\bc\in \R^{k-\ell, \ell}$, $\Xi(x)\bc$ is a hypersurface in $\R^{k+1}$ with flat normal bundle, $x$ is a line of curvature coordinate system, and the first fundamental form $\I= \sum_{i=1}^k g_{ii} \rd x_i^2$ satisfies the condition that $\sum_{i=1}^n \e_i g_{ii}$ is equal to the length of $\bc$ in $\R^{k-\ell, \ell}$. In particular, if $\bc$ is a null vector in $\R^{k-\ell, \ell}$ then $\Xi\bc$ is an isothermic$_\ell$ hypersurface  (as defined in \ref{fg}).
\item For any $\bc_1, \bc_2\in \R^{k-\ell, \ell}$, $\Xi(x)\bc_1\mapsto \Xi(x)\bc_2$ is a Combescure O-transform.  
\item Let $C=(\bc_1, \ldots, \bc_k)$ be a constant matrix in $GL(k)$,  and
$$Y=(Y_1, \ldots, Y_k)= \Xi C, \quad \hat Y=(\hat Y_1, \ldots, \hat Y_k)=\hat \Xi C.$$  
Then:
\ben
\item[(a)] $Y, \hat Y$ are $k$-tuples in $\R^{k+1}$ of type $\R^{k-\ell, \ell}$,
\item[(b)] If all columns of $C$ are null vectors in $\R^{k-\ell, \ell}$, then $Y_i, \hat Y_i$ are isothermic$_\ell$ hypersurfaces in $\R^{k+1}$ for $1\leq i\leq k$. 
\item[(c)] $E(x,0), \hat E(x,0)\in O(k+1)\times O(k-\ell, \ell)$.
\item[(d)] Write 
$ E(\cdot,0)=\bpm g_1&0\\ 0 & g_2\epm$, $ \hat E(\cdot,0)= \bpm \hat g_1 & 0\\ 0 & \hat g_2\epm$, 
$$g_1= (e_1, \ldots, e_{k+1}), \quad \hat g_1= (\hat e_1, \ldots, \hat e_{k+1}), $$
and $\hat W= (q_1, \ldots, q_{k+1})^t$. 
Then 
$$\hat Y_i= Y_i -\frac{\hat Z^t Jg_2^{-1}\bc_i}{s} \sum_{j=1}^{ k+1} q_j e_j.$$ 
\item[(e)] $Y_i(x)\mapsto \hat Y_i(x)$ is a Ribaucour transform for $1\leq i\leq k$. In fact, we have
$$\hat Y_i - \frac{\hat Z^tg_2^{-1}\bc_i}{sq_j} \hat e_j = Y_i - \frac{\hat Z^tg_2^{-1}\bc_i}{s q_j} e_j$$
for all $1\leq i, j\leq k+1$.
\een
\een
\ethm

Note that a solution of the $\frac{O(n+k -\ell, \ell)}{O(n)\times O(k-\ell, \ell)}$-system gives rise to a family of isothermic$_\ell$ $k$-submanifolds in $\R^n$ parametrized by the null cone of $\R^{k-\ell, \ell}$ and any two submanifolds in this family are related by Combescure O-transforms.  But for the converse, we need to have $k$ or $k-1$  
isothermic$_\ell$ $k$-submanifolds in $\R^n$ related by Combescure O-transforms to construct a solution of the $\frac{O(n+k -\ell, \ell)}{O(n)\times O(k-\ell, \ell)}$-system. This is because Theorem \ref{eg} (1)-(3) hold for any Combescure O-map $Y=(Y_1, \ldots, Y_m):\R^k\to \calM_{n\times m}$. So the connection $o(k-\ell, \ell)$-valued  $1$-form $\tau=\d F^t- JF\d J$ has $m$ parallel sections, and $\tau$ is flat if $m=k-1$ or $m=k$.

\end{document}